\documentclass[12pt]{article}

\usepackage{latexsym}

\begin{document}


\newcommand{\deu}[1]{}
\newcommand{\eng}[1]{#1}

\newcommand{\dontdothat}[1]{}

%

\newcommand{\EE}{\mathop{\rm I\! E}\nolimits}
\newcommand{\E}{\mathop{\rm E}\nolimits}
\newcommand{\I}{\mathop{\rm Im  }\nolimits}
\newcommand{\R}{\mathop{\rm Re  }\nolimits}
\newcommand{\CC}{\mathop{\rm C\!\!\! I}\nolimits}
\newcommand{\FF}{\mathop{\rm I\! F}\nolimits}
\newcommand{\KK}{\mathop{\rm I\! K}\nolimits}
\newcommand{\LL}{\mathop{\rm I\! L}\nolimits}
\newcommand{\MM}{\mathop{\rm I\! M}\nolimits}
\newcommand{\NN}{\mathop{\rm I\! N}\nolimits}
\newcommand{\PP}{\mathop{\rm I\! P}\nolimits}
\newcommand{\QQ}{\mathop{\rm I\! Q}\nolimits}
\newcommand{\RR}{\mathop{\rm I\! R}\nolimits}
\newcommand{\ZZ}{\mathop{\rm Z\!\!Z}\nolimits}
\newcommand{\integer}{\mathop{\rm int}\nolimits}
\newcommand{\erf}{\mathop{\rm erf}\nolimits}
\newcommand{\diag}{\mathop{\rm diag}\nolimits}
\newcommand{\fl}{\mathop{\rm fl}\nolimits}
\newcommand{\eps}{\mathop{\rm eps}\nolimits}
\newcommand{\var}{\mathop{\rm var}\nolimits}


\newcommand{\myownlab}[1]{{\label{#1}}}

\newcommand{\ring}{{\cal K}}

\newcommand{\pfeil}{\rightarrow}

\newcommand{\HA}{{\rm HA}}
\newcommand{\HB}{{\AK^{k}}}
\newcommand{\BA}{{\rm BA}}
\newcommand{\TB}{{\rm TB}}	   \newcommand{\TA}{{\rm TA}}
\newcommand{\ZA}{{\rm ZA}}
\newcommand{\ZB}{{\rm ZB}}
\newcommand{\HBB}{{\rm\cal AK}}
\newcommand{\BAB}{{\rm\cal BA}}
\newcommand{\BBfull}{{\rm BB}}
\newcommand{\AKBB}{{\rm BB^{k}}}
\newcommand{\AK}{{\rm AK}}
\newcommand{\BB}{{\rm BB^{\ast}}}
\newcommand{\BBB}{{\rm\cal BB^{\ast}}}
\newcommand{\AKBBB}{{\rm\cal BB}^{(k)}}
\newcommand{\AKBBs}{{\rm {\widetilde{BB^{(k)}}}}}
\newcommand{\BBs}{{\rm {\widetilde{BB^{\ast}}}}}
\newcommand{\tr}{{\rm tr}}
\newcommand{\Tr}{{\rm Tr}}
\newcommand{\iso}{\stackrel{\sim}{=}}

\newcommand{\kat}{{\cal C}}
\newcommand{\rmat}{{\cal R}}
\newcommand{\oalg}{{\cal A}}
\newcommand{\falg}{{\cal F}}
\newcommand{\eich}{{\cal G}}
\newcommand{\hilb}{{\cal H}}
\newcommand{\calm}{{\cal M}}
\newcommand{\mod}{{\cal M}}

\newcommand{\bigrho}{\rho_\oplus}
\newcommand{\bigphi}{\phi_\oplus}
\newcommand{\qquer}{{\overline{q}}}

\newenvironment{bew}{Proof:}{\hfill$\Box$}
\newtheorem{bem}{Remark}
\newtheorem{bsp}{Example}
\newtheorem{axiom}{Axiom}
\newtheorem{satz}{Proposition}
\newtheorem{de}[satz]{Definition}
\newtheorem{lede}[satz]{Definition and Lemma}
\newtheorem{lemma}[satz]{Lemma}
\newtheorem{kor}[satz]{Corollary}
\newtheorem{theo}[satz]{Theorem}
\newtheorem{hypo}[satz]{Hypothesis}

\newcommand{\sbegin}[1]{\small\begin{#1}}
\newcommand{\send}[1]{\end{#1}\normalsize}

\sloppy

\newcommand{\lang}{1}

\title{An Ariki-Koike like extension of the \\ Birman-Murakami-Wenzl Algebra}
\author{Reinhard H\"aring-Oldenburg\\ 
          Mathematisches Institut\\ Bunsenstr. 3-5\\
 37073 G\"ottingen, Germany\\
email: haering@cfgauss.uni-math.gwdg.de}
\date{20.8.97}
\maketitle

\begin{abstract} We introduce an Ariki-Koike like 
extension of the Birman-Murakami-Wenzl Algebra and show it
to be semi-simple.  This algebra  supports a faithful Markov trace that
gives rise to link invariants of closures of Coxeter type B braids.
\end{abstract}

\section{Introduction}
The theory of quantum invariants of links nowaday rests on a broad
theory that includes quantum groups, their centraliser algebras
and tensor categories. It is the ultimate goal of the 'Knot Theory and
Root Systems' programme initiated in \cite{tD1}	to carry over this
theory to the braid groups associated to the other root systems.
The greatest progress  sofar has been taken for the braid group of
Coxeter type B where  the notions of quasi triangular
Hopf algebra and monoidal categories have been defined and 
nontrivial examples have been found \cite{DHO}, \cite{DHO2}, 
\cite{rhobcat}. Furthermore, Temperley-Lieb algebras 
and Hecke algebras have been studied intensively for this root system.
In the present paper we continue the study of generalisations
of the Birman-Murakami-Wenzl algebra \cite{we2}, \cite{rhobmw}.

Every Coxeter diagram defines a braid group that is an infinite
covering of its Coxeter group. 
The braid group $\ZB_n$ of Coxeter type B has 
generators $\tau_i,i=0,1,\ldots n-1$.
Generators $\tau_i,i\geq1$  satisfy the relations of 
Artin's braid group (which is
the braid group of Coxeter type A): 
$\tau_i\tau_j=\tau_j\tau_i$ if $|i-j|>1$, and
$\tau_i\tau_j\tau_i=\tau_j\tau_i\tau_j$ if $|i-j|=1$. 
The generator $\tau_0$ has relations
\begin{eqnarray}
\tau_0\tau_1\tau_0\tau_1&=&\tau_1\tau_0\tau_1\tau_0\label{vierzopf}\\
\tau_0\tau_i&=&\tau_i\tau_0\quad i\geq2
\end{eqnarray} 

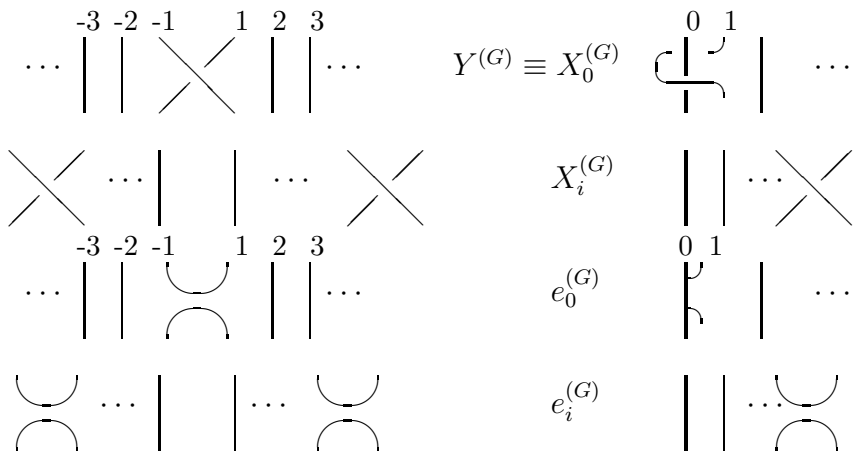
\begin{figure}[ht]
\begin{picture}(140,60)
\put(59,50){\mbox{$Y^{(G)}\equiv X_0^{(G)}$}}

\linethickness{0.2mm}
\put(2,50){\mbox{$\cdots$}}
\put(9,56){\mbox{{\small -3}}}
\put(14,56){\mbox{{\small -2}}}
\put(19,56){\mbox{{\small -1}}}
\put(40,56){\mbox{{\small 3}}}
\put(35,56){\mbox{{\small 2}}}
\put(30,56){\mbox{{\small 1}}}

\put(10,45){\line(0,1){10}}
\put(15,45){\line(0,1){10}}
\put(20,55){\line(1,-1){10}}
\put(20,45){\line(1,1){4}}
\put(26,51){\line(1,1){4}}
\put(35,45){\line(0,1){10}}
\put(40,45){\line(0,1){10}}
\put(42,50){\mbox{$\cdots$}}

\linethickness{0.4mm}
\put(90,45){\line(0,1){3}}
\put(90,50){\line(0,1){5}}
\linethickness{0.2mm}
\put(88,51){\oval(4,4)[l]}
\put(88,49){\line(1,0){5}}
\put(93,47){\oval(4,4)[tr]}
\put(93,55){\oval(4,4)[br]}
\put(95,56){\mbox{{\small 1}}}
\put(90,56){\mbox{{\small 0}}}

\put(100,45){\line(0,1){10}}

\put(107,50){\mbox{$\cdots$}}

\put(72,35){\mbox{$X_i^{(G)}$}}

\put(0,40){\line(1,-1){10}}
\put(0,30){\line(1,1){4}}
\put(6,36){\line(1,1){4}}
\put(45,40){\line(1,-1){10}}
\put(45,30){\line(1,1){4}}
\put(51,36){\line(1,1){4}}
\put(20,30){\line(0,1){10}}
\put(30,30){\line(0,1){10}}
\put(13,35){\mbox{$\cdots$}}
\put(35,35){\mbox{$\cdots$}}
\linethickness{0.4mm}
\put(90,30){\line(0,1){10}}
\linethickness{0.2mm}
\put(95,30){\line(0,1){10}}
\put(98,35){\mbox{$\cdots$}}
\put(102,40){\line(1,-1){10}}
\put(102,30){\line(1,1){4}}
\put(108,36){\line(1,1){4}}

\put(72,20){\mbox{$e_0^{(G)}$}}

\linethickness{0.2mm}
\put(2,20){\mbox{$\cdots$}}
\put(9,26){\mbox{{\small -3}}}
\put(14,26){\mbox{{\small -2}}}
\put(19,26){\mbox{{\small -1}}}
\put(40,26){\mbox{{\small 3}}}
\put(35,26){\mbox{{\small 2}}}
\put(30,26){\mbox{{\small 1}}}

\put(10,15){\line(0,1){10}}
\put(15,15){\line(0,1){10}}
\put(25,15){\oval(8,8)[t]}
\put(25,25){\oval(8,8)[b]}
\put(35,15){\line(0,1){10}}
\put(40,15){\line(0,1){10}}
\put(42,20){\mbox{$\cdots$}}

\linethickness{0.4mm}
\put(90,15){\line(0,1){10}}
\linethickness{0.2mm}
\put(90,17){\oval(4,4)[tr]}
\put(90,25){\oval(4,4)[br]}
\put(93,26){\mbox{{\small 1}}}
\put(89,26){\mbox{{\small 0}}}

\put(100,15){\line(0,1){10}}

\put(107,20){\mbox{$\cdots$}}

\put(72,5){\mbox{$e_i^{(G)}$}}
\linethickness{0.4mm}
\put(90,0){\line(0,1){10}}
\linethickness{0.2mm}
\put(95,0){\line(0,1){10}}
\put(98,5){\mbox{$\cdots$}}
\put(106,0){\oval(8,8)[t]}
\put(106,10){\oval(8,8)[b]}
\put(20,0){\line(0,1){10}}
\put(30,0){\line(0,1){10}}
\put(12,5){\mbox{$\cdots$}}
\put(32,5){\mbox{$\cdots$}}
\put(5,0){\oval(8,8)[t]}
\put(5,10){\oval(8,8)[b]}
\put(45,0){\oval(8,8)[t]}
\put(45,10){\oval(8,8)[b]}

\end{picture}
\caption{\label{generat} The graphical interpretation of the 
generators as symmetric tangles (on the left)
and as cylinder tangles (on the right)}
\end{figure}

The braid group   $\ZB_n$ may be 
graphically interpreted (cf. figure \ref{generat}) as 
symmetric braids or cylinder braids \cite{tD3}: The symmetric picture
shows it as the group  of 
braids with $2n$ strands (numbered $-n,\ldots,-1,1,\ldots,n$) which are
fixed under a 180 degree rotation about the middle axis. 
In the cylinder picture one adds a single fixed line (indexed $0$)
on the left and
obtains $\ZB_n$ as the group of braids with $n$ strands that may 
surround this fixed line. 
The generators $\tau_i,i\geq0$ are mapped 
to the diagrams $X^{(G)}_i$ given in figure \ref{generat}.
More generally, tangles of B-type may be defined. The special case
of tangles without crossings is the B-type 
 Temperley-Lieb algebra $\TB_n$ that has been introduced
by tom Dieck in \cite{tD1}.

The Ariki-Koike Algebra is the quotient of the group algebra of $\ZB_n$ 
where the images $X_i$ of the generators $\tau_i$ for $i\geq1$
fulfil quadratic relations while $X_0$ satisfies a polynomial of
arbitrary degree. The Hecke Algebra of
B type is a special case where  $X_0$ satisfies also a quadratic relation.

The standard Birman-Murakami-Wenzl algebra $\BA_n$
of type A imposes cubic relations
on its generators in a way that enables its interpretation  as an algebra
of tangles with a skein relation that comes from the Kauffman polynomial.

Thus it is natural to define an Ariki-Koike like extension 
of the BMW algebra $\AKBB_n$ that  
contains a generator $Y$ as image of $\tau_0$ that satisfies 
$\prod_{i=0}^{k-1}(Y-p_i)=0$. 
The special case $k=2$ has been called restricted B-type BMW algebra
and has been studied in \cite{rhobmw}.

The current interest in the study of B type braid groups has several
origins. Closing B type braids yields links that can be interpreted as links 
in a solid torus \cite{lamb}
and Markov traces on group algebras of $\ZB_n$ hence 
allow the calculations of invariants of such links (cf. end of section 
\ref{mainsec}). 
Braid groups of all finite root systems further act as symmetries
on the corresponding quantum groups \cite{lu}.	The B braid group
occurs furthermore in several physical situations \cite{rhoref},
\cite{rhopott}.	The general idea  is that
the B type braids allow to treat with  knot 
theoretic methods also  physical models with a boundary. 
The $\tau_0$ generator is interpreted as a reflection at the boundary.

We now outline the structure of the paper and point out the main results.
After a short review of the Birman-Wenzl algebra of A-type 
we go on to define the 
Ariki-Koike-Birman-Murakami-Wenzl algebra of B-type 
$\AKBB_n$ in section 	\ref{AKdefsec} and list a  number of 
fundamental relations. They are used extensively in section \ref{wordsec} 
to determine a partial normal form 	of words in $\AKBB_n$.
Section \ref{connections} shows how to obtain the 
Ariki-Koike algebra as a quotient of $\AKBB_n$.
Furthermore, it investigates the B type Temperley Lieb sub-algebra. 

Section \ref{clsec}  introduces the graphical 
interpretation of our algebra
and studies its classical limit. 
The construction of a Markov trace fills section \ref{trsec}. 

The main theorem of this paper is contained in section \ref{mainsec}.
We prove that $\AKBB_n$ is semi-simple in the generic case and show how its 
simple components can be enumerated in terms of Young diagrams.
The Bratteli diagram is given and we show that the Markov trace is faithful.
As an application a generalisation of the Kauffman polynomial to 
links in the solid torus is discussed.

  \unitlength1mm

  \unitlength1mm
 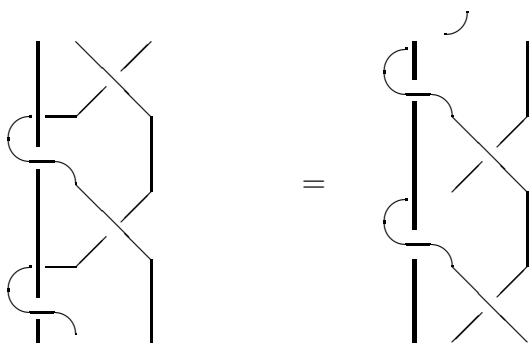
\begin{figure}[ht]
\begin{picture}(150,44)

\linethickness{0.4mm}
\put(10,40){\line(0,-1){14}}
\put(10,23){\line(0,-1){17}}
\put(10,3){\line(0,-1){3}}
\linethickness{0.2mm}

\put(15,40){\line(1,-1){10}}
\put(25,40){\line(-1,-1){4}}
\put(15,30){\line(1,1){4}}
\put(15,30){\line(-1,0){4}}

\put(9,27){\oval(6,6)[l]}
\put(9,24){\line(1,0){3}}
\put(12,21){\oval(6,6)[tr]}
\put(25,30){\line(0,-1){10}}

\put(15,21){\line(1,-1){10}}
\put(15,10){\line(1,1){4}}
\put(25,20){\line(-1,-1){4}}

\put(25,11){\line(0,-1){11}}
\put(15,10){\line(-1,0){4}}
\put(9,7){\oval(6,6)[l]}
\put(9,4){\line(1,0){3}}
\put(12,1){\oval(6,6)[tr]}


\put(45,20){\mbox{$=$}}

\linethickness{0.4mm}
\put(60,40){\line(0,-1){5}}
\put(60,32){\line(0,-1){17}}
\put(60,11){\line(0,-1){11}}
\linethickness{0.2mm}

\put(59,36){\oval(6,6)[l]}
\put(59,33){\line(1,0){3}}
\put(62,30){\oval(6,6)[tr]}
\put(75,40){\line(0,-1){10}}
   \put(64,44){\oval(6,6)[br]}

\put(65,30){\line(1,-1){10}}
\put(65,20){\line(1,1){4}}
\put(75,30){\line(-1,-1){4}}

\put(59,16){\oval(6,6)[l]}
\put(59,13){\line(1,0){3}}
\put(62,10){\oval(6,6)[tr]}
\put(75,20){\line(0,-1){10}}

\put(65,10){\line(1,-1){10}}
\put(65,0){\line(1,1){4}}
\put(75,10){\line(-1,-1){4}}

\end{picture}

\caption{\label{vierzopfbild} 
Relation (\ref{vierzopf}) in the cylinder picture}
\end{figure}

{\em Algebraic preliminaries}:

We collect some simple results from algebra that will be needed later on.

Our first topic is the specialisation of the ground ring
of an algebra.
Let  $R$ and $R'$ be rings and let $\pi:R\rightarrow R'$ 
denote an epimorphism. Let $S$ be a set that is considered as a set
of generators. The set of words over this alphabet is denoted by 
$S^\ast$. The free $R$ module with basis $S^\ast$ is denoted by
$\langle S\rangle_R$. We then have an induced mapping
\[
\pi_\ast:\langle S\rangle_R\rightarrow
\langle S\rangle_{R'}, 
\qquad \sum_ir_iw_i\mapsto\sum_i\pi(r_i)w_i,
\quad r_i\in R, w_i\in S^\ast.\] 
An alternative description of $\pi_\ast$ is possible in
terms of tensor products by viewing  $R'$ to be a $R$ module
(with the action given by $\pi$). We then have
$\langle S\rangle_R\otimes_RR'=\langle S\rangle_{R'}$.
Now, let
$M=\{\sum_jm_{ij}w_j\mid m_{i,j}\in R,w_j\in S^\ast\}$
be a set of relations and $\overline{M}$ the ideal generated 
by these elements in	 $\langle S\rangle_R$. 
The projection is denoted by $p:\langle S\rangle_R
\rightarrow\langle S\rangle_R/\overline{M},a\mapsto
a+\overline{M}$. Furthermore, let $\overline{M}'$ be the
ideal generated by $\pi_\ast(M)$
and denote by $p'$ the corresponding projection. 
It  follows  that
\[\pi_\ast p=p'\pi_\ast\]

We now turn to results about the construction of integral domains.
\begin{satz} 
Let  $K$ be an infinite field. An ideal $I\subset K[x_1,\ldots,x_n]$
is prime if, and only if its affine variety $V(I)$  is irreducible.

$V(I)$ is irreducible if it can be rationally parametrised 
\[x_i=\frac{f_i(t_1,\ldots,t_m)}{g_i(t_1,\ldots,t_m)},
\qquad i=1,\ldots n\]
Here  $f_i,g_i$ denote polynomials.
\end{satz}

\begin{kor}		 \label{primsatz}
Let  $K$ be an infinite field. Assume the ideal
$I\subset K[x_1,\ldots,x_n]$
is generated by  $m$ polynomials $h_i$ which may  
in the field of fractions 
$K(x_1,\ldots,x_n)$ be solved uniquely for
$x_1,\ldots,x_m$:
\[
x_i=\frac{f_i(x_{m+1},\ldots,x_n)}{g_i(x_{m+1},\ldots,x_n)}
\qquad i=1,\ldots,m	\]
Then $I$ is prime.
\end{kor}

\begin{satz}\label{esymsatz}
Let $K$ denote an algebraically closed field and let 
$\sigma_i,i=1,\ldots,n$ be the elementary symmetric polynomials
in $K[x_1,\ldots,x_n]$. Then the  system of equations 
$\sigma_i=y_i$ over 
$K[x_1,\ldots,x_n,y_1,\ldots,y_n]$ can be solved for the $x_i$.
\end{satz}

\section[Definition of the Ariki-Koike-BMW-Algebra]{The Definition 
of the Ariki-Koike-Birman-Mu\-ra\-ka\-mi-Wenzl-Algebra}	
  \myownlab{AKdefsec}\label{defsec}

This section introduces a generalisation of the 
Birman-Murakami-Wenzl that is related to the B-type braid group.
Because the algebras of Ariki and Koike appear as quotients 
we call our algebra an Ariki-Koike-BMW algebra.
We set off by recalling the definition of the ordinary BMW algebra.

\begin{de}	  \myownlab{defba}
Let  $R$ denote an integral domain with units $x,\lambda\in R$
such that with a further element $\delta\in R$  the relation 
$(1-x)\delta=\lambda-\lambda^{-1}$ holds.
The Birman-Murakami-Wenzl (BMW) algebra
$\BA_n(R)$ is generated by 
$X_1,\ldots,X_{n-1},e_1,\ldots,e_{n-1}$ and relations:
\begin{eqnarray}
X_iX_j&=&X_jX_i\qquad |i-j|>1\myownlab{def2}\\
X_iX_jX_i&=&X_jX_iX_j\qquad |i-j|=1\myownlab{def3}\\
X_ie_i&=&e_iX_i=\lambda e_i\myownlab{def4}\\
e_iX_j^{\pm1}e_i&=&\lambda^{\mp1}e_i\qquad |i-j|=1\myownlab{def5}\\
e_i^2&=&xe_i\myownlab{lem1a}\\
X_i^{-1}&=&X_i-\delta+\delta e_i\myownlab{lem1d}\\
e_ie_j&=&e_je_i\qquad |i-j|>1\myownlab{lem1f}\\
e_iX_jX_i&=&X_j^\pm X_i^\pm e_j\qquad |i-j|=1\myownlab{lem1h}\\
e_ie_je_i&=&e_i\qquad |i-j|=1\myownlab{lem1l}
\end{eqnarray}
\end{de}

\begin{lemma}\label{balemma1}
\begin{eqnarray}
X_i^2&=&1+\delta X_i-\delta\lambda e_i\myownlab{lem1e}\\
X_i^3&=&X_i^2(\lambda+\delta)+X_i(1-\lambda\delta)-\lambda\myownlab{lem1b}\\
X_i^{-2}&=&1+\delta^2-\delta X_i+\delta(\lambda^{-1}-\delta)e_i
  = 1-\delta X_i^{-1}+\delta\lambda^{-1} e_i\myownlab{lem1q}\\
X_i^{-1}X_j^{\pm1}X_i&=&X_jX_i^{\pm1}X_j^{-1}\qquad |i-j|=1\myownlab{lem1g}\\
X_i^{\pm1}e_je_i&=&X_j^{\mp1}e_i\qquad |i-j|=1\myownlab{lem1k}\\
e_ie_jX_i^{\pm1}&=&e_iX_j^{\mp1}\qquad |i-j|=1\myownlab{lem1kk}\\
e_iX^\pm_jX^\pm_i&=&e_ie_j\qquad |i-j|=1\myownlab{lem1m}\\
X_i^\pm X_j^\pm e_i&=&e_je_i\qquad |i-j|=1\myownlab{lem1mm}\\
X_ie_jX_i^{-1}&=&X_j^{-1}e_iX_j\qquad |i-j|=1\myownlab{lem1i}\\
X_ie_jX_i&=&X_j^{-1}e_iX_j^{-1}\qquad |i-j|=1\myownlab{lem1p}
\end{eqnarray}
\end{lemma}
\begin{bew}
(\ref{lem1e})-(\ref{lem1q}) are simple restatements of  (\ref{lem1d}).

(\ref{lem1g}): $X_iX_jX_i=X_jX_iX_j\Rightarrow
                X_jX_iX_j^{-1}=X_i^{-1}X_jX_i\Rightarrow
                X_iX_j^{-1}X_i^{-1}=X_j^{-1}X_i^{-1}X_j$

(\ref{lem1k}): $X_i^\pm e_je_i=X_j^\mp X_j^\pm X_i^\pm e_je_i
\stackrel{(\ref{lem1h})}{=}
X_j^\mp e_iX_j^\pm X_i^\pm e_i=\lambda^\pm X_j^\mp e_iX_j^\pm e_i
\stackrel{(\ref{def5})}{=}
X_j^\mp e_i$

(\ref{lem1kk}): $e_ie_jX_i^\pm=e_ie_jX_i^\pm X_j^\pm X_j^\mp=
e_iX_i^\pm X_j^\pm e_i X_j^\mp=
\lambda^\pm e_iX_j^\pm e_iX_j^\mp=e_iX_j^\mp $ 

(\ref{lem1m}): Using (\ref{lem1l}), (\ref{lem1k}) and 
         (\ref{lem1h}) we have
\begin{eqnarray*}
e_iX_j^\pm X_i^\pm &=&e_ie_je_iX^\pm _jX^\pm _i=
e_iX_i^{\mp}X_i^\pm e_je_iX^\pm _jX_i^\pm \\
&=&e_iX_i^{\mp}X_j^{\mp}e_iX^\pm _jX^\pm _i=
e_iX_i^{\mp}X_j^{\mp}X^\pm _jX^\pm _ie_j=e_ie_j
\end{eqnarray*}
 
(\ref{lem1mm}) is shown similarly to (\ref{lem1m}). 
(\ref{lem1i},\ref{lem1p}) follow from 
(\ref{lem1h}).

\end{bew}

\begin{lemma} \label{balemma}
If $\delta$ is a unit in  $R$ the algebra $\BA(R)$ 
is isomorphic to the algebra generated by invertible 
$X_1,\ldots,X_{n-1}$ and relations (\ref{def2})-(\ref{def5}).
The element $e_i$ is now defined by
\begin{eqnarray}
e_i&:=&1-\frac{X_i-X_i^{-1}}{\delta}\myownlab{edef}\qquad i=1,\ldots,n-1
\end{eqnarray}
\end{lemma}
\begin{bew}
(\ref{lem1a}): $e_i^2=(1-\delta^{-1}(X_i-X_i^{-1}))e_i=
  e_i-\delta^{-1}(\lambda e_i-\lambda^{-1}e_i)=xe_i $

(\ref{lem1f}): follows from (\ref{def2}) using (\ref{edef})

(\ref{lem1h}): from (\ref{edef}) and (\ref{lem1g})

(\ref{lem1l}): The middle  $e_j$ is replaced by (\ref{edef}):
\[e_ie_je_i=xe_i-\delta^{-1}(e_iX_je_i-e_iX_j^{-1}e_i)=
    (1-\delta^{-1}(\lambda-\lambda^{-1}))e_i-
            \delta^{-1}(\lambda^{-1}e_i-\lambda e_i)=e_i\]

\end{bew}

We now define  our generalised algebra.

\begin{de} \label{definAK}
Fix $k\in\NN$ and let $x,\lambda,\kappa,p_0,\ldots,p_{k-1}\in R$  be units
and let $\delta,A_1,\ldots,A_{k-1}\in R$	be some further elements.
Assume that the relation $(1-x)\delta=\lambda-\lambda^{-1}$ holds. 
 The Ariki-Koike-BWM-Algebra
 on $n$ strands  $\AKBB_n(R)$ is defined as
 $R$ algebra generated by
$Y,X_1,\ldots,X_{n-1},e_1,\ldots,e_{n-1}$ and the relations of
the  Birman-Murakami-Wenzl-Algebra $\BA_n$ and
\begin{eqnarray}
X_1YX_1Y&=&YX_1YX_1\myownlab{AKdef6}\\
YX_i&=&X_iY\qquad i>1\myownlab{AKdef9}\\
YX_1Ye_1&=&\kappa e_1\myownlab{AKdef8}\\
0&=&\prod_{i=0}^{k-1}(Y-p_i)\myownlab{AKdef7} \\
e_1Y^ie_1&=&A_ie_1\qquad i=1,\ldots,k-1\myownlab{AKdef10}
\end{eqnarray}
\end{de}

Relation (\ref{AKdef10}) suggests to define $A_0:=x$.

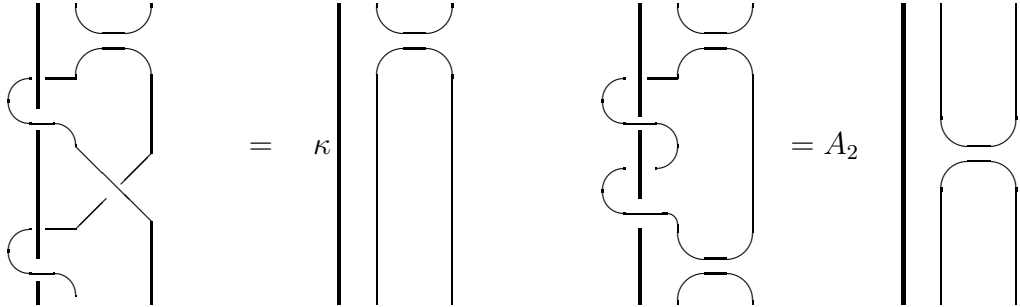
\begin{figure}[ht]
\begin{picture}(150,44)

\linethickness{0.4mm}
\put(10,40){\line(0,-1){14}}
\put(10,23){\line(0,-1){17}}
\put(10,3){\line(0,-1){3}}
\linethickness{0.2mm}

\put(20,30){\oval(10,8)[t]}
\put(20,40){\oval(10,8)[b]}

\put(15,30){\line(-1,0){4}}

\put(9,27){\oval(6,6)[l]}
\put(9,24){\line(1,0){3}}
\put(12,21){\oval(6,6)[tr]}
\put(25,30){\line(0,-1){10}}

\put(15,21){\line(1,-1){10}}
\put(15,10){\line(1,1){4}}
\put(25,20){\line(-1,-1){4}}

\put(25,11){\line(0,-1){11}}
\put(15,10){\line(-1,0){4}}
\put(9,7){\oval(6,6)[l]}
\put(9,4){\line(1,0){3}}
\put(12,1){\oval(6,6)[tr]}


\put(38,20){\mbox{$=\quad\kappa$}}

\linethickness{0.4mm}
\put(50,40){\line(0,-1){40}}
\linethickness{0.2mm}

\put(60,30){\oval(10,8)[t]}
\put(60,40){\oval(10,8)[b]}

\put(55,30){\line(0,-1){30}}
\put(65,30){\line(0,-1){30}}


\linethickness{0.4mm}
\put(90,40){\line(0,-1){15}}
\put(90,23){\line(0,-1){9}}
\put(90,0){\line(0,1){10}}
\linethickness{0.2mm}

\put(100,30){\oval(10,8)[t]}
\put(100,40){\oval(10,8)[b]}

\put(105,30){\line(0,-1){20}}

\put(91,30){\line(1,0){4}}
\put(88,27){\oval(6,6)[l]}
\put(88,15){\oval(6,6)[l]}
\put(92,21){\oval(6,6)[r]}
\put(88,24){\line(1,0){4}}
\put(88,12){\line(1,0){5}}
\put(93,10){\oval(4,4)[tr]}

\put(100,0){\oval(10,8)[t]}
\put(100,10){\oval(10,8)[b]}

\put(110,20){\mbox{$= A_2$}}

\linethickness{0.4mm}
\put(125,40){\line(0,-1){40}}
\linethickness{0.2mm}

\put(135,15){\oval(10,8)[t]}
\put(135,25){\oval(10,8)[b]}

\put(130,25){\line(0,1){15}}
\put(140,25){\line(0,1){15}}

\put(130,15){\line(0,-1){15}}
\put(140,15){\line(0,-1){15}}

\end{picture}

\caption{\label{minzopfbild} 
Relation (\ref{AKdef8}) (on the left) and relation (\ref{AKdef10})
(on the right) in the cylinder picture}
\end{figure}

These relations are motivated by our intended graphical
interpretation.
Section \ref{clsec} will give precise definitions of the
graphical version of the algebra. Here we only shed some
light on the interpretation of the relations. 
(\ref{AKdef6}) is the four braid relation (\ref{vierzopf})
which is visualised in figure \ref{vierzopfbild}.
Relation (\ref{AKdef9}) stems from the braid group as well.
Relation (\ref{AKdef8}) is visualised in figure \ref{minzopfbild}. 
The graphical calculus suggests to take either
 $\kappa=1$ or $\kappa=\lambda$ (depending on the
 precise ribbon graph which $Y$ should represent).
Disconnected components of a graph may be eliminated
using (\ref{AKdef10}). Finally, (\ref{AKdef7}) is motivated
 by algebraic considerations.

  \unitlength1mm

The generic ground ring for our algebra is 
a quotient of a Laurent polynomial ring. 
We denote by  $R[x]$ the polynomial ring and by 
$R\{x\}$  the Laurent ring in $x$ over $R$.
\begin{eqnarray}
\label{r0def}
R_0:=\CC[\delta,A_1,\ldots,A_{k-1}]
\{p_0,\ldots,p_{k-1},x,\kappa,\lambda\}/
(x\delta-\delta-\lambda^{-1}+\lambda)
\end{eqnarray}
The ring's dependence on $k$ is not written explicitly.

\begin{bem} \myownlab{AKinvolution} 
There is an involution of $\AKBB_n(R_0)$ as a $\CC$ algebra
such that 
\begin{equation}
X_i^\ast:=X_i^{-1},e_i^\ast:=e_i,Y^\ast:=Y^{-1},
\delta^\ast:=-\delta,x^\ast:=x,\lambda^\ast:=\lambda^{-1}, 
p_i^\ast:=p_i^{-1},\kappa^\ast:=\kappa^{-1}
\nonumber
\end{equation}
$A_i^\ast$ has to be calculated from $e_1Y^{-i}e_1=A_i^\ast e_1$
using the following formulas.
\end{bem}

\begin{de}
Let $q_{k-1},\ldots,q_0$ be the signed elementary symmetric 
polynomials in   $p_0,\ldots,p_{k-1}$ such that:
\begin{equation}
Y^k=\sum_{i=0}^{k-1}q_iY^i	   \myownlab{AKykexp}
\end{equation}
\end{de}

Note that $q_0=(-1)^{k-1}\prod_ip_i$ is invertible.
We calculate $Y^{-1}$:
\begin{eqnarray}
Y^{-1}&=&\sum_{i=0}^{k-1}\qquer_iY^i\quad
\mbox{ with }\quad \qquer_{k-1}=q_0^{-1}\quad\qquer_{i-1}=-q_iq_0^{-1}
\end{eqnarray}
The coefficients are determined uniquely if the $Y^i$ 
are linearly independent.

Iterating one obtains expressions $\overline{Q}_{i,j}$ such that:
\begin{equation}
Y^{-i}=\sum_{j=0}^{k-1}\overline{Q}_{i,j}Y^j
\end{equation}
Acting with the  involution $\ast$ one obtains
\begin{equation}
Y^{i}=\sum_{j=0}^{k-1}Q_{i,j}Y^{-j}\qquad
\mbox{with}\quad Q_{i,j}=\overline{Q}^\ast_{i,j}
\end{equation}

The following definitions will prove useful later on.
\begin{eqnarray}
Y_i&:=&X_{i-1}X_{i-2}\cdots X_1YX_1^{-1}\cdots X_{i-2}^{-1}X_{i-1}^{-1}\\
Y^{(m)}_i&:=&X_{i-1}X_{i-2}\cdots X_1Y^mX_1\cdots X_{i-2}X_{i-1}\\
Y'_i&:=&Y^{(1)}_i=X_{i-1}X_{i-2}\cdots X_1YX_1\cdots X_{i-2}X_{i-1}
\end{eqnarray}

The next lemma collects a stock of relations that
show among other things that the most important properties of $Y$ can
be shifted to other strands.
\begin{lemma}\myownlab{lemma2}
\begin{eqnarray}
Y_i^k&=&\sum_{j=0}^{k-1}q_jY_i^j\myownlab{AKlem2aa}\\
Y_i^{-1}&=&\sum_{j=0}^{k-1}\qquer_jY^j_i\myownlab{AKlem2aaa}\\
0&=&[X_1YX_1Y,\{Y,e_1,X_1\}]\myownlab{AKlem2b}\\
0&=&[Y_i,X_j]=[Y_i,e_j]\qquad j\neq i,i-1\myownlab{AKlem2h}\\
0&=&[Y^{(m)}_i,X_j]=[Y^{(m)}_i,e_j]\qquad j\neq i,i-1\myownlab{AKlem2hh}\\
Y^{(m)}_{i+1}X_i^{-1}&=&X_iY^{(m)}_i \myownlab{AKyl1}\\
Y_{i+1}X_i&=&X_iY_i\myownlab{AKyl2}\\
X_iY_iX_iY_i&=&Y_iX_iY_iX_i\myownlab{AKlem2zopf}\\
X_iY'_iX_iY'_i&=&Y'_iX_iY'_iX_i\myownlab{AKlem2zopfs}\\
\kappa e_i&=&e_iY_iX_iY_i=Y_iX_iY_ie_i\myownlab{AKlem2c}\\
\kappa e_i&=&e_iY'_iX_iY'_i=Y'_iX_iY'_ie_i
\myownlab{AKlem2cs}\\
e_iY_i^me_i&=&A_me_i\myownlab{AKlem2ds}\\
Y'_iY'_j&=&Y'_jY'_i\myownlab{AKlem2f}\\
Y_iY_{i-1}^{-1}&=&
  Y^{-1}_{i-1}X_{i-1}^{-1}Y_{i-1}X_{i-1}\myownlab{AKlem2yyinv}\\
Y_ie_{i-1}&=&\kappa \lambda^{-1}Y_{i-1}^{-1}e_{i-1}  \myownlab{AKlem2q1}\\
e_{i-1}Y_i&=&\kappa \lambda e_{i-1}Y_{i-1}^{-1}-\lambda\delta e_{i-1}Y_{i-1}+
  \delta\lambda A_1e_1  \myownlab{AKlem2q2}\\
e_{i-1}Y'_i&=&\kappa\lambda e_{i-1}{Y'_{i-1}}^{-1}\myownlab{AKlem2m}\\
Y'_ie_{i-1}&=&\kappa\lambda {Y'_{i-1}}^{-1}e_{i-1}\myownlab{AKlem2mm}\\
X_iY_{i+1}&=&  X_iY_i-\delta Y_i+\delta Y_ie_i+\delta Y_{i+1}-
\kappa\delta\lambda e_iY_i^{-1}+\nonumber\\&&\delta^2\lambda e_iY_i-
\delta^2\lambda A_1e_i\myownlab{AKlem2q3}\\  
Y_{i+1}^lX_i&=&X_iY_i^l\myownlab{AKlem2q4}\\
e_iY_i^lX_i&=&\kappa\lambda e_iY_i^{l-1}X_iY_i^{-1}-
 \kappa\delta\lambda e_iY_i^{l-2}+
  \kappa\delta\lambda A_{l-1}e_iY_i^{-1} \myownlab{AKlem2eyx}\\
X_iY_i^le_i&=&\kappa\lambda Y_i^{-1}X_iY_i^{l-1}e_i-
  \delta\kappa\lambda Y_i^{l-2}e_i+\kappa\delta\lambda A_{l-1}Y_i^{-1}e_i
    \myownlab{AKlem2xye}
\end{eqnarray}
\end{lemma}
\begin{bew}

(\ref{AKlem2b}): Using (\ref{AKdef6}), one has
$X_1X_1YX_1Y=X_1YX_1YX_1$ and hence $X_1YX_1Y$ commutes with 
$X_1$. Thus it also commutes with  $X_1^{-1}$ and $e_1$. 

(\ref{AKlem2h},\ref{AKlem2hh}): For $j\geq i+1$ 
commutativity follows from (\ref{def2},\ref{AKdef9})
and for $j\leq i-1$ 
it is an application of (\ref{def3}).

(\ref{AKyl1}) and (\ref{AKyl2}) are trivial.
(\ref{AKlem2zopf}), (\ref{AKlem2zopfs})
are shown by induction. The induction step for (\ref{AKlem2zopf}) 
reads:
\begin{eqnarray*}
Y_iX_iY_iX_i&=&X_{i-1}Y_{i-1}X_{i-1}^{-1}X_iX_{i-1}Y_{i-1}X_{i-1}^{-1}X_i\\
&=&X_{i-1}Y_{i-1}X_iX_{i-1}X_{i}^{-1}Y_{i-1}X_{i-1}^{-1}X_i\\
&=&X_{i-1}X_iY_{i-1}X_{i-1}Y_{i-1}X_i^{-1}X_{i-1}^{-1}X_i\\
&=&X_{i-1}X_iY_{i-1}X_{i-1}Y_{i-1}X_{i-1}X_i^{-1}X_{i-1}^{-1}\\
&=&X_{i-1}X_iX_{i-1}Y_{i-1}X_{i-1}Y_{i-1}X_i^{-1}X_{i-1}^{-1}\\
&=&X_{i}X_{i-1}X_{i}Y_{i-1}X_{i-1}X_i^{-1}Y_{i-1}X_{i-1}^{-1}\\
&=&X_{i}X_{i-1}Y_{i-1}X_{i}X_{i-1}X_i^{-1}Y_{i-1}X_{i-1}^{-1}\\
&=&X_{i}X_{i-1}Y_{i-1}X_{i-1}^{-1}X_{i}X_{i-1}Y_{i-1}X_{i-1}^{-1}\\
&=&X_{i}Y_{i}X_{i}Y_{i}
\end{eqnarray*}
The induction step for (\ref{AKlem2zopfs}) is almost identical.

(\ref{AKlem2c},\ref{AKlem2cs}): The inductive proofs start
from (\ref{AKdef8}) and its mirror version:
\[\lambda e_1YX_1Y=e_1X_1YX_1Y\stackrel{(\ref{AKlem2b})}{=}
X_1YX_1Ye_1=\lambda YX_1Ye_1=\kappa\lambda e_1\]

In the induction step for (\ref{AKlem2c}) 
equation (\ref{lem1h}) is used to eliminate $e_{i+1}$ in terms of  $e_i$:
\begin{eqnarray*}
Y^{(1)}_{i+1}X_{i+1}Y^{(1)}_{i+1}e_{i+1}
&=&X_iY^{(1)}_iX_iX_{i+1}X_iY^{(1)}_iX_iX_i^{-1}X_{i+1}^{-1}e_iX_{i+1}X_i\\
&=&X_iY^{(1)}_iX_{i+1}X_iX_{i+1}Y^{(1)}_iX_{i+1}^{-1}e_iX_{i+1}X_i\\
&=&X_iX_{i+1}Y^{(1)}_iX_iX_{i+1}X_{i+1}^{-1}Y^{(1)}_ie_iX_{i+1}X_i\\
&=&X_iX_{i+1}Y^{(1)}_iX_iY^{(1)}_ie_iX_{i+1}X_i=
\kappa X_iX_{i+1}e_iX_{i+1}X_i=\kappa e_{i+1}
\end{eqnarray*}
The induction step for (\ref{AKlem2cs}) is:
\begin{eqnarray*}
e_{i+1}Y_{i+1}X_{i+1}Y_{i+1}
&=&e_{i+1}X_iY_iX_i^{-1}X_{i+1}X_iY_iX_i^{-1}
=e_{i+1}X_iY_iX_{i+1}X_iX_{i+1}^{-1}Y_iX_i^{-1}\\
&=&e_{i+1}X_iX_{i+1}Y_iX_iY_iX_{i+1}^{-1}X_i^{-1}\
=X_iX_{i+1}e_iY_iX_iY_iX_{i+1}^{-1}X_i^{-1}\\
&=&\kappa X_iX_{i+1}e_iX_{i+1}^{-1}X_i^{-1}=\kappa e_{i+1}\\
\end{eqnarray*}

(\ref{AKlem2ds}): The proof is by induction.
\begin{eqnarray*}
e_iY_i^me_i&=&e_iX_{i-1}Y^m_{i-1}X_{i-1}^{-1}e_i=
e_ie_{i-1}X_i^{-1}Y^m_{i-1}X_{i-1}^{-1}e_i\\
&=&e_ie_{i-1}Y^m_{i-1}X_i^{-1}X_{i-1}^{-1}e_i
=e_ie_{i-1}Y^m_{i-1}e_{i-1}X_i^{-1}X_{i-1}^{-1}\\
&=&A_me_ie_{i-1}X_i^{-1}X_{i-1}^{-1}\stackrel{(\ref{lem1kk})}{=}A_me_i
\end{eqnarray*}

(\ref{AKlem2f}): $[Y,Y^{(1)}_1]=[Y,Y^{(1)}_2]=0$ is trivial. For
$i>1$ the claim follows by induction: 
$[Y,Y^{(1)}_i]=0\Rightarrow[Y,Y^{(1)}_{i+1}]=[Y,X_iY^{(1)}_iX_i]=0$. 
In the general case  $[Y^{(1)}_j,Y^{(1)}_i]$ we may assume   $j<i$. 
Using (\ref{AKlem2h}) the induction step is: 
$[Y^{(1)}_j,Y^{(1)}_i]=[X_{j-1}Y^{(1)}_{j-1}X_{j-1},Y^{(1)}_i]=0$.

(\ref{AKlem2yyinv}) is a consequence of (\ref{AKlem2zopf}).

(\ref{AKlem2q1}):
\begin{eqnarray*}
 Y_ie_{i-1}&=&X_{i-1}Y_{i-1}X_{i-1}^{-1}e_{i-1}=
   \lambda^{-1}X_{i-1}Y_{i-1}e_{i-1}
 =\kappa\lambda^{-1}Y_{i-1}^{-1}e_{i-1}
\end{eqnarray*}

(\ref{AKlem2q2}):
\begin{eqnarray*}
e_{i-1}Y_i&=&e_{i-1}X_{i-1}Y_{i-1}X_{i-1}^{-1}=
\lambda	e_{i-1}Y_{i-1}X_{i-1}^{-1}\\
&=&
\lambda e_{i-1}Y_{i-1}X_{i-1}-\delta\lambda e_{i-1}Y_{i-1}+
\delta\lambda e_{i-1}Y_{i-1}e_{i-1}	\\
&=&\kappa\lambda e_{i-1}Y_{i-1}^{-1}-\lambda\delta e_{i-1}Y_{i-1}+
  \delta\lambda A_1e_1 
  \end{eqnarray*}

(\ref{AKlem2m},\ref{AKlem2mm}) are shown in the following way:
\[e_{i-1}Y^{(1)}_i=e_{i-1}X_{i-1}Y^{(1)}_{i-1}X_{i-1}=
\lambda e_{i-1}Y^{(1)}_{i-1}X_{i-1}Y^{(1)}_{i-1}Y_{i-1}^{(1)-1}=
\kappa\lambda e_{i-1}Y_{i-1}^{(1)-1}\]

(\ref{AKlem2q3}): $X_iY_{i+1}=X_i^2Y_iX_i^{-1}=
Y_iX_i^{-1}+\delta Y_{i+1}-\delta\lambda e_iY_iX_i^{-1}=
X_iY_i-\delta Y_i+\delta Y_ie_i+\delta Y_{i+1}-
\kappa\delta\lambda e_iY_i^{-1}+\delta^2\lambda e_iY_i-
\delta^2\lambda A_1e_i$

(\ref{AKlem2q4}) is trivial. We show (\ref{AKlem2eyx}):
\begin{eqnarray*}
e_iY_i^lX_i&=&e_iY_i^lX_iY_iX_iX_i^{-1}Y_i^{-1}=
e_iX_iY_iX_iY_i^lX_i^{-1}Y_i^{-1}\\
&=&\lambda e_iY_iX_iY_i^lX_i^{-1}Y_i^{-1}
=\kappa\lambda e_iY_i^{l-1}X_i^{-1}Y_i^{-1}  \\
&=&\kappa\lambda e_iY_i^{l-1}X_iY_i^{-1}-
\kappa\delta\lambda e_iY_i^{l-1}Y_i^{-1}+
\kappa\delta\lambda e_iY_i^{l-1}e_iY_i^{-1}\\
&=&\kappa\lambda e_iY_i^{l-1}X_iY_i^{-1}-\kappa\delta\lambda e_iY_i^{l-2}+
\kappa\delta\lambda A_{l-1}e_iY_i^{-1}
\end{eqnarray*}
\end{bew}

\begin{bem} \myownlab{reliinvol} 
 $X_i^\dagger:=X_{n-i}, Y^\dagger:=Y_n$ defines an involution $\dagger$
 on  $\AKBB_n$.
\end{bem}
\begin{bew}
All relations that depend only on one index or on the
absolute difference of two indices are obviously compatible.
We check  (\ref{AKdef10}):
\begin{eqnarray*}
(e_1Y^ie_1-A_ie_1)^\dagger&=&e_{n-1}Y_n^ie_{n-1}-A_ie_{n-1}\\
&=&
e_{n-1}X_{n-1}Y_{n-1}^iX_{n-1}^{-1}e_{n-1}-A_ie_{n-1}=0
\end{eqnarray*}
Relation (\ref{AKdef8}) is preserved as well: 
\begin{eqnarray*}
(YX_1Ye_1-\kappa e_1)^\dagger&=&
e_{n-1}Y_nX_{n-1}Y_n-\kappa e_{n-1}\\
&=&
e_{n-1}X_{n-1}Y_{n-1}X_{n-1}^{-1}X_{n-1}X_{n-1}Y_{n-1}X_{n-1}^{-1}-
\kappa e_{n-1}\\&=&
\lambda e_{n-1}Y_{n-1}X_{n-1}Y_{n-1}X_{n-1}^{-1}-\kappa e_{n-1}\\
&=&
\lambda \kappa e_{n-1}X_{n-1}^{-1}-\kappa e_{n-1}=0
\end{eqnarray*}
\end{bew}

\begin{bem}	 \myownlab{topdowninvol}
The relations show that there is a further
involution  $a\mapsto\overline{a}$ which fixes all generators.
\end{bem}

\section{Relations to other Knot algebras}\myownlab{connections}

The  $e_i$ together with a projector $e_0$ on the $p_0$ 
eigenvalue of $Y$ generate a  sub-algebra that is a homomorphic image of
a Type-B-Temperley-Lieb algebra. 
The quotient by the ideal generated by this sub-algebra is 
isomorphic to the  Ariki-Koike algebra.
For specific parameter values one may also
obtain the A-type BMW algebra as a quotient.

\begin{lemma}\myownlab{jdef}
Let $J_n$ be the ideal generated by  $Y_n-p_0$.
Every other $Y_i-p_0,i=1,\ldots,n$ 
generates the same ideal and the quotient 
$R/(\kappa-\lambda p_0^2,xp_0^i-A_i)\otimes_R\AKBB_n(R)/J_n$ 
is isomorphic to the A-type BMW algebra
$\BA_n(R/(\kappa-\lambda p_0^2,xp_0^i-A_i))$.
\end{lemma}
\begin{bew}
The first claim is a consequence of the definition of the $Y_i$.
The specialisation of  $\kappa$ and  $A_i$ is necessary
since in the quotient one obtains
$0=e_1YX_1Y-\kappa e_1=e_1p_0X_1p_0-\kappa e_1=e_1(p_0^2\lambda-\kappa)$
and $A_ie_1=e_1Y^ie_1=p_0^ie_1$.
The remaining relations present no further restrictions.
\end{bew}

\begin{de}
$I_n$ denotes the ideal generated by $e_{n-1}$ in $\AKBB_n$.
\end{de}

As we shall see, the quotient by this ideal is an Ariki-Koike algebra.
\begin{de}				\myownlab{heckedef}
 $\AK^{k}_n$ denotes the  Ariki-Koi\-ke algebra \cite{ariki}
with generators
$X_0$, $X_1$, $\ldots,$ $X_{n-1}$ and parameters 
$\delta,p_i,i=0,\ldots k-1$ and relations:
\begin{eqnarray}
X_0X_1X_0X_1&=&X_1X_0X_1X_0\myownlab{hdef1}\nonumber\\
X_iX_j&=&X_jX_i\qquad |i-j|>1\myownlab{hdef2}\nonumber\\
X_iX_jX_i&=&X_jX_iX_j\qquad |i-j|=1\myownlab{hdef3}\nonumber\\
X_i^2&=& \delta X_i+1\quad i\geq 0\myownlab{hdef4}\nonumber\\
0&=&\prod_{i=0}^{k-1}(X_0-p_i)	  \myownlab{hdef5}	\nonumber
\end{eqnarray}
\end{de}
We use a slightly different normalisation of the parameters 
than Ariki and Koike did. From their work we need the result
that $\AK^{k}_n$ is semi-simple.
The proof of the following lemma is now trivial.
\begin{lemma}					\myownlab{heckelemma}
 $I_n$ is generated by any of the $e_i$ and the 
 quotient by it is isomorphic to $AK^{k}_n$.
\end{lemma}

Of some interest in knot theoretical applications is the
projector on the eigenvalue   $p_0$ of $Y$. 
Such a projector is given by $\prod_{i=1}^{k-1}(Y-p_i)$,
but we prefer to work with sums of $Y^i$.
\begin{de} 	Let $e_0=\sum_{i=0}^{k-1}\alpha_i Y^i$ 
be a projector on the eigenvalue $p_0$ of $Y$, i.e. it
satisfies 
 $e_0Y=Ye_0=p_0e_0$.
\end{de}

\begin{lemma}\begin{eqnarray}
\alpha_{k-1}&=&\alpha_0p_0q_0^{-1}\qquad 
  \alpha_{j-1}=p_0\alpha_j-\alpha_{k-1}q_j\myownlab{alpharec}\\
e_0^2&=&x_0e_0\qquad x_0:=\sum_{i=0}^{k-1}\alpha_ip_0^i\\
e_1e_0e_1&=&x'_0e_1\qquad x'_0:=\sum_{i=0}^{k-1}\alpha_i A_i 
\end{eqnarray}\end{lemma}
The proofs are simple.

The modified B-Temperley-Lieb Algebra 
(\cite{tD1}, \cite{rhopott}) $\TB'_n$ is defined by generators
  $e_0,e_1,\ldots,e_{n-1}$,
parameters $c,c',d$ and relations $e_0^2=ce_0,e_i^2=de_i,
e_je_l=e_le_j,e_ie_je_i=e_i,e_1e_0e_1=c'e_1,1\leq i,j \leq n-1,
0\leq l \leq n-1, |i-j|=1, |j-l|>1$.
Obviously, we have a morphism $\TB'_n\rightarrow\AKBB_n$.

\section{Two strands $n=2$ and ground rings}

The algebra $\AKBB_n(R)$ is in general not semisimple.
This section studies conditions that suffice to make
 $B(R):=\AKBB_2(R)$ semisimple.
For the sake of notaional convenience we omit the index $1$ of
$e_1$ and $X_1$.

The parameters of the algebra cannot be choosen
independendly. Note for example that both
 $e=\kappa e YXY$ and $Y^k=\sum_iq_iY^i$ fix the length
of  $Y$.  
\begin{de}				 \myownlab{defconsist}
Define a ring  $R_1$ as a quotient   
$R_1:=R_0/c$ of $R_0$. The ideal $c\subset R_0$ is generated by 
 $k$ Laurent polynomials  that are obtained by the following
 procedure: Expand $Ye-\kappa X_1^{-1}Y^{-1}e$ using
(\ref{AKdef10}), (\ref{lem1d}),
(\ref{AKykexp}) and  (\ref{AKlem2xye}) into a linear combination 
$\sum_{i=0}^{k-1} h_iY^ie$. 
The coefficient  $h_i$ of this sum are the generators of
 $c=(h_0,\ldots,h_{k-1})$.
\end{de}

\begin{lemma}
\label{negexp}
The expansion of the expression in the definition of $c$
in terms of  $Y^{-i}e$ defines polynomials
 $h'_i$ such that
$Ye-\kappa X_1^{-1}Y^{-1}e=\sum_{i=0}^{k-1} h'_iY^{-i}e$.
They generate the same ideal: $c=(h'_0,\ldots,h'_{k-1})$.
Furthermore, $c$ is closed under the involution
$c^\ast=c$. 
\end{lemma}
\begin{bew}
We have  $h_j =\sum_i\overline{Q}_{i,j}h'_i$ and 
        $h'_j=\sum_iQ_{i,j}h_i$. This implies equality of
		both the ideals generated by these sets of polynomials.
		The second claim follows from $h'_i=h_i^\ast$.
\end{bew}

To shed some light on the ideal  $c$ we first note that
 (\ref{AKlem2xye}) implies:
\[
X_1Y^me=\lambda^{m-1}\kappa^mY^{-m}e+\sum_{s=1}^{m-1}
\kappa^s\lambda^s\delta(A_{m-s}Y^{-s}e-Y^{m-2s}e)
\]
This renders the defining equation into the form
\begin{eqnarray}	 \myownlab{consisteq}
Ye+\delta\kappa Y^{-1}e-
\kappa\delta\sum_m\overline{q}_mA_me_1\myownlab{cosisteq}\\
=\kappa\sum_m\overline{q}_m  
\left(\lambda^{m-1}\kappa^mY^{-m}e+\sum_{s=1}^{m-1}
\kappa^s\lambda^s\delta(A_{m-s}Y^{-s}e-Y^{m-2s}e)  \right)\nonumber
\end{eqnarray}

We now introduce a ring that will become relevant later on as the ring of the
classical limit of the algebra. At this stage we need it
purely as a tool.

\begin{de}
 The ideal
$J_c\subset R_1$ is given by
$J_c:=(\kappa-1,\lambda-1,q-1,q_0-1,q_1,\ldots,q_{k-1})$.
Set $R_c:=R_1/J_c$.
\end{de}

According to proposition \ref{esymsatz} the equation for the  $q_i$
are solvable. Hence the ring 
$R_c$ is nontrivial. 
The same polynomials  $(\kappa-1,\lambda-1,q-1,q_0-1,q_1,\ldots,q_{k-1})$
define an ideal in  $R_0$. It contains $c$
since after dividing by $J_c$ we have $Y^{-1}=Y^{k-1},
\overline{q}_{k-1}=1,\overline{q}_i=0$
and hence (\ref{consisteq})  becomes trivial.	
It follows that $R_c$ is the quotient of $R_0$ by $J_c$.

The ring $R_1$ plays a special role in the conctruction of a 
$B$-module.
Using  (\ref{AKdef10}),(\ref{AKykexp}) and (\ref{AKlem2xye})
we see that the ideal $I_2$ is spanned $Y^ieY^j,i,j=0,1,\ldots,k-1$.
\begin{de} Let $R$ be as in  definition \ref{definAK}.	
Let $V:=V(R)$ be the free $R$-module of dimension  $k$. 
The basis is denoted by
 $b_i,0\leq i<k$. 
 $V$ is turned into a module 
 of the free algebra generated by  $e,X,Y$ by the following 
 definitions:
\begin{eqnarray*}
e.b_i&:=&A_ib_0\qquad \mbox{ with } A_0=x\\
Y.b_{k-1}&:=&\sum_{j=0}^{k-1}q_jb_j\qquad Y.b_i:=b_{i+1}\\
X.b_0&:=&\lambda b_0\qquad X.b_1:=\kappa Y^{-1}.b_0\\
X.b_i&:=&\kappa\lambda(Y^{-1}.X.b_{i-1}-\delta b_{i-2}+
\delta A_{i-1}Y^{-1}.b_0)\quad i\geq2
\end{eqnarray*}
\end{de}

The definition of this  action is guided by the desire that 
it should factor over $B(R)$.

 $Y^{-1}$ and $X^{-1}$  shall act by their expansions
 in terms of  $Y^i$
(implying  $(Y^{-1}).=(Y.)^{-1}$ ),
resp. $X,e,1$.  
It turns out, however, that $V$ is not in general a
 $B$-module.
 Most relations are easy to check but two of them may not hold:
 (a) $XYXY=YXYX$ and (b) $X^2=1+\delta X-\delta\lambda e$.
Relation (b) is equivalent to  $(X^{-1}).=(X.)^{-1}$.

\begin{lemma} $V(R_1)$ is a $B(R_1)$-module.
\end{lemma}
\begin{bew}
For the ring  $R_1$ one has by its construction:
 \begin{equation}\label{b1prop}
 b_1=Y.b_0=\kappa X^{-1}.Y^{-1}.b_0
 \end{equation}
On   $b_0$ relation (b) holds trivially. We check (a):
\begin{eqnarray*}
X.Y.X.Y.b_0&=&  X.Y.X.b_1=\kappa X.Y.Y^{-1}.b_0
=\kappa X.b_0=\kappa\lambda b_0\\
&=&\lambda\kappa Y.Y^{-1}.b_0=\lambda Y.X.b_1=Y.X.Y.X.b_0
\end{eqnarray*}
Furthermore, we check the inverse of  (a):
\begin{eqnarray*}
X^{-1}.Y^{-1}.X^{-1}.Y^{-1}.b_0
&\stackrel{(\ref{b1prop})}{=}&
\kappa^{-1}X^{-1}.Y^{-1}.Y.b_0=\kappa^{-1}X^{-1}.b_0=
\kappa^{-1}\lambda^{-1} b_0\\
Y^{-1}.X^{-1}.Y^{-1}.X^{-1}.b_0&=&
\lambda^{-1}Y^{-1}.X^{-1}.Y^{-1}.b_0\stackrel{(\ref{b1prop})}{=}
\kappa^{-1}\lambda^{-1}Y^{-1}.Y.b_0=\lambda^{-1}\kappa^{-1}b_0
\end{eqnarray*}   
(\ref{b1prop}) enables us to write  for all $i=0,\ldots,k-1$:
\begin{equation}\label{sternrel}
X.Y.b_{i-1}=X.b_i=\kappa\lambda Y^{-1}.X^{-1}.b_{i-1}
\end{equation}
Here we used the convention that $b_{-1}=Y^{-1}.b_0$.
The case $i=0$ follows from  (\ref{b1prop}), 
the case  $i=1$ is trivial, and the cases $i>1$
are simple rewritings of  the action of $X$. 

Now we can start the inductive proof that
(a) and (b) hold on all basis vectors. 
The induction assumption 
 $H_i$ is: Relations (a) and (b) hold on $b_{i-1}$. 
 We show that the inverse of relation (a) holds on $b_{i-1}$. 
\begin{eqnarray*}
Y^{-1}.X^{-1}.Y^{-1}.X^{-1}.b_{i-1}&\stackrel{(\ref{sternrel})}{=}&
\kappa^{-1}\lambda^{-1}X.Y.\kappa^{-1}\lambda^{-1}X.Y.b_{i-1}\\
&=&\kappa^{-2}\lambda^{-2}X.Y.X.Y.b_{i-1}
\stackrel{H_i}{=}\kappa^{-1}\lambda^{-1}b_{i-1}\\
X^{-1}.Y^{-1}.X^{-1}.Y^{-1}.b_{i-1}
&=&X^{-1}.Y^{-1}.X^{-1}.b_{i-2}
\stackrel{(\ref{sternrel})}{=}
\kappa^{-1}\lambda^{-1}X^{-1}.X.b_{i-1}\\
&\stackrel{H_i}{=}&
\kappa^{-1}\lambda^{-1}b_{i-1}
\end{eqnarray*}
We now check  (b):
\begin{eqnarray*}
X^{-1}.X.b_i&\stackrel{(\ref{sternrel})}{=}&
\kappa\lambda X^{-1}.Y^{-1}.X^{-1}.b_{i-1}
=\kappa\lambda Y.Y^{-1}.X^{-1}.Y^{-1}.X^{-1}.b_{i-1}\\
&=&\kappa\lambda Y.X^{-1}.Y^{-1}.X^{-1}.Y^{-1}.b_{i-1}
\stackrel{(\ref{sternrel})}{=}Y.X^{-1}.X.b_{i-1}
\stackrel{H_i}{=}b_i
\end{eqnarray*}
Finally, we look at (a):
\begin{eqnarray*}
Y.X.Y.X.b_i&=&\kappa\lambda Y.X.Y.Y^{-1}.X^{-1}.b_{i-1}
=\kappa\lambda Y.X.X^{-1}.b_{i-1}\stackrel{H_i}{=}
\kappa\lambda Y.b_{i-1}=\kappa\lambda b_i\\
X.Y.X.Y.b_i&\stackrel{(\ref{sternrel})}{=}&
\kappa\lambda X.Y.Y^{-1}.X^{-1}.b_i=
\kappa\lambda X.X^{-1}.b_i=\kappa\lambda b_i
\end{eqnarray*}
\end{bew}

\begin{de}
$U_m:= {\rm span}_{R_1}\{Y^ieY^m\mid i=0,\ldots,k-1\}$
\end{de}
\begin{lemma}
Each $U_m$ is a $B(R_1)$-module isomorphic to  $V=V(R_1)$.
They have pairwise trivial intersections.
\end{lemma}
\begin{bew}
We show that the map $\varrho:V(R_1)\rightarrow B(R_1), b_i\mapsto Y^ie$
defines a module isomorphism of $V$ and $U_0$.
It is a surjection of $R_1$-modules, and, by the above lemma,
a morphism of $B(R_1)$-modules.
It remains to check injectivity.
Suppose we had $0=\sum_i\alpha_iY^ie,\alpha_i\in R_1$.
Applying this to $b_0$ we obtain 
$0=x\sum_i\alpha_ib_i$. Now,  $x$ is invertible, and hence
all the $\alpha_i$ have to vanish.
Thus we have shown that ${\rm span}_{R_1}\{Y^ie_1\}$ 
is a free $R_1$ module. The same is true for the isomorphic 
$B(R_1)$ modules $U_m$.
Now, we are going to show that the $e_1$ ideal  
${\rm span}\{Y^ie_1Y^j\}$ as a whole is a free $R_1$ module.
It suffices to show that the $U_m$ form a direct sum decomposition,
i.e. that $m\neq r\Rightarrow
U_m\cap U_r=\{0\}$.
Since  $Y$ is invertible, it suffices to show 
for $m\geq1$ that $U_m\cap U_0=\{0\}$.
Assume there is a non zero element  
$a$ in the intersection of  $U_0$ and  $U_m$:
\[a=\sum_ia_iY^ie_1=\sum_ib_iY^ie_1Y^m\]
Multiplying from the right with  $e_1$ the righthand side is 
mapped to $U_0$ and we may compare the coefficients in its basis:
$a_i=x^{-1}A_mb_i$. Thus
\[a=x^{-1}A_m\sum_ib_iY^ie_1=\sum_ib_iY^ie_1Y^m\]
Hence  $x^{-1}A_maY^{-m}=aY^{-m}Y^m$, that is $aY^{-m}$ is an eigenvector
of $Y^m$ to the eigenvalue  $x^{-1}A_m$. Now, $Y^m$ and $x$ are invertible
and hence  $A_m$ is a unit in $R_1$ and furthermore in any
quotient of $R_1$. However, we have already defined the quotient $R_c$ 
in which $A_m$ is obviously not invertible.  
\end{bew}

\begin{lemma}
$R_1$ is an integral domain.
\end{lemma}
\begin{bew}
We have to show that the ideal  $c$ is prime.
Due to proposition \ref{primsatz} it suffices to show that
the defining equations may be solved uniquely in the field of fractions.

Using lemma \ref{negexp}
we consider the coefficients of $Y^{-i}e_1, i=1,\ldots,k-2$.
They form a triangular (and hence a soluble) system 
of equations in the variables $A_i$. To solve it one first determines
$A_{1}$ from the coefficient of $Y^{-(k-2)}e_1$. Secondly,  $A_2$ 
is calculated from the coefficient of $Y^{-(k-3)}e_1$.
We end with $A_{k-2}$ and  $Y^{-1}e_1$.
Thereafter the coefficient of  $e_1$ can be used to isolate
$A_{k-1}$ which appears just once in this expression.
It remains to investigate the coefficient of $Y^{-(k-1)}e_1$.
It is $q_0=\kappa q_0^{-1}\lambda^{k-2}\kappa^{k-1}
-\kappa\sum_{m=0}^{k-1}\sum_{s=1}^{m-1}\overline{q}_m\kappa^s\lambda^s
(\lambda-\lambda^{-1})(1-x)^{-1}
{\rm Coeff}(Y^{m-2s}e_1,Y^{-(k-1)}e_1)$. 
This can be solved for  $x$.
\end{bew}

Note that the proof of this lemma breaks down if one chooses
to specify $\kappa=\lambda^{-1}$ since then we can't be sure
that the coefficient of    $(1-x)^{-1}$ is non-zero. 

\begin{de}
Let $K_1$ denote the field of fractions of $R_1$.
\end{de}

The quotient of  $B(R_1)$ by the ideal $I_2$ 
is isomorphic to the Ariki-Koike algebra $\AK^{k}_2(R_1)$
which is a free module over any integral domain
\cite{ariki}.
We summarise:
\begin{lede}\myownlab{yemodul}
$B(R_1)$ is a free  $R_1$ module and the subset
$Y^ie_1Y^j$ is linearly independent.
\end{lede}

\section{The word problem in $\AKBB_n$}	\myownlab{wordsec}

In this section we single out a set of words in standard form
that linearly generate $\AKBB_n$. However, this does not lead to a
linear basis of $\AKBB_n$ but it is fundamental to the following 
analysis.

\begin{satz}\myownlab{wortsatz}
Every element in $\AKBB_n$ is a linear combination of words of the form
$w_1\gamma w_2$ where $w_i\in\AKBB_{n-1}$ and
$\gamma\in\Gamma_n:=\{1,e_{n-1},X_{n-1},Y_{n}^j,j=1,\ldots,k-1\}$
The same is true if in $\Gamma_n$ the generators $X_{n-1}$ or 
$Y_n$ or both are  replaced by their inverses.
\end{satz}
\begin{bew} We prove the proposition by induction.
The case $n=1$ is trivial and $n=2$ can also be verified easily.

Let $w_0\gamma_0w_1\gamma_1\cdots w_m\gamma_mw_{m+1}\in\AKBB_n,
w_i\in\AKBB_{n-1}$ be an arbitrary word.
It suffices to show that any two neighbouring $\gamma_i$ 
can be combined together.
Hence the situation we have to investigate is $w=\gamma_1w_1\gamma_2,
w_1\in\AKBB_{n-1},\gamma_1,\gamma_2\in\Gamma_n$. 
By induction hypothesis we have
$w_1=u_1\alpha u_2,u_i\in\AKBB_{n-2},\alpha\in\Gamma_{n-1}$ 
and hence
$w=\gamma_1u_1\alpha u_2\gamma_2=u_1\gamma_1\alpha\gamma_2 u_2$. 
Thus it suffices to investigate
$w'=\gamma_1\alpha\gamma_2$. The cases $\gamma_1=1$ or $\gamma_2=1$
are trivial. 
We now investigate in turn the four possible values of $\alpha$.

1. Case $\alpha=1$: 
The following table gives the relation that allows to reduce the product
$\gamma_1\gamma_2$ to the standard form of the proposition.
\[\begin{array}{c||c|c|c}
\gamma_1\backslash\gamma_2 & Y_n^j & e_{n-1} & X_{n-1}\\\hline\hline
Y_n^l &(\ref{AKlem2aa}) & (\ref{P2a1}) & (\ref{AKlem2q4})\\\hline
e_{n-1} &(\ref{P2a2}) & (\ref{lem1a})&(\ref{def4})\\\hline
X_{n-1} & (\ref{P2a3})& (\ref{def4}) & (\ref{lem1e})
\end{array}\]

\begin{equation}\myownlab{P2a1}
Y_n^le_{n-1}=X_{n-1}Y_{n-1}^lX_{n-1}^{-1}e_{n-1}=
\lambda^{-1}X_{n-1}Y_{n-1}^le_{n-1}
\mbox{ apply (\ref{AKlem2xye}) recursively}
\end{equation}

\begin{eqnarray}\myownlab{P2a2}
e_{n-1}Y_n^j&=&\lambda e_{n-1}Y_{n-1}^jX_{n-1}^{-1}\\
&=&\lambda e_{n-1}Y_{n-1}^jX_{n-1}-\delta\lambda e_{n-1}Y_{n-1}^j
+\delta\lambda A_j e_{n-1}
\nonumber
\end{eqnarray}
The first term is reduced by applying (\ref{AKlem2eyx})
recursively.

\begin{eqnarray}\myownlab{P2a3}
X_{n-1}Y_n^j&=&X_{n-1}^2Y_{n-1}^jX_{n-1}^{-1}\\
&=&Y_{n-1}^jX_{n-1}^{-1}-
\delta\lambda e_{n-1} Y_{n-1}^j X_{n-1}^{-1}
+\delta X_{n-1}Y_{n-1}^jX_{n-1}^{-1}\nonumber \\
&=&Y_{n-1}^jX_{n-1}-\delta Y_{n-1}^j+\delta Y_{n-1}^je_{n-1}\\&&
-\delta\lambda 
(e_{n-1} Y_{n-1}^j X_{n-1}-\delta e_{n-1} Y_{n-1}^j
+\delta A_j e_{n-1} )
+\delta Y_{n}^j
\end{eqnarray}
Again, one needs (\ref{AKlem2eyx}) for recursive reduction.

2. Case $\alpha=X_{n-2}$: 
\[\begin{array}{c||c|c|c}
\gamma_1\backslash\gamma_2 & Y_n^j & e_{n-1} & X_{n-1}\\\hline\hline
Y_n^l & =X_{n-2}Y_n^{j+l} \quad (\ref{AKlem2aa})&
  =X_{n-2}Y_n^le_{n-1}\quad (\ref{P2a1})
   &=X_{n-2}Y_n^lX_{n-1}\quad (\ref{AKlem2q4})\\\hline
e_{n-1} &=e_{n-1}Y_n^jX_{n-2}\quad(\ref{P2a2}) &
  (\ref{def5})&(\ref{lem1m})\\\hline
X_{n-1} & =X_{n-1}Y_n^jX_{n-2}\quad(\ref{P2a3})& (\ref{lem1mm}) & (\ref{def3})
\end{array}\]

3. Case $\alpha=e_{n-2}$: 
\[\begin{array}{c||c|c|c}
\gamma_1\backslash\gamma_2 & Y_n^j & e_{n-1} & X_{n-1}\\\hline\hline
Y_n^l & =e_{n-2}Y_n^{l+j}\quad (\ref{AKlem2aa})& 
=e_{n-2}Y_n^le_{n-1}\quad (\ref{P2a1})& 
=e_{n-2}Y_n^lX_{n-1}\quad (\ref{AKlem2q4})  \\\hline
e_{n-1} & =e_{n-1}Y_n^je_{n-2}\quad(\ref{P2a2}) 
  & (\ref{lem1l})&(\ref{lem1kk})\\\hline
X_{n-1} &  X_{n-1}Y_n^je_{n-2}\quad(\ref{P2a3})  
  & (\ref{lem1k}) & (\ref{lem1p})
\end{array}\]

4. Case $\alpha=Y_{n-1}^m$:
\[\begin{array}{c||c|c|c}
\gamma_1\backslash\gamma_2 & Y_n^j & e_{n-1} & X_{n-1}\\\hline\hline
Y_n^l & \ast  &  \mbox{ like } (\ref{P2a5}) & \ast \\\hline
e_{n-1} & (\ref{P2a5}) & (\ref{AKlem2ds})&  (\ref{AKlem2eyx}) \\\hline
X_{n-1} & \ast &(\ref{AKlem2xye})  & (\ref{flcc}) 
\end{array}\]

\begin{eqnarray}\myownlab{P2a5}
e_{n-1}Y_{n-1}^mY_n^j&=&e_{n-1}Y_{n-1}^mX_{n-1}Y_{n-1}^jX_{n-1}^{-1}\\
&\stackrel{(\ref{AKlem2eyx})}{\in}& 
{\rm span}\{ e_{n-1}Y_{n-1}^s\mid 0\leq s<k\}Y_{n-1}^jX_{n-1}^{-1}
\nonumber\\
 &\stackrel{(\ref{AKlem2eyx})}{\subseteq}& 
 {\rm span}\{ e_{n-1}Y_{n-1}^s\mid 0\leq s<k\}
\end{eqnarray}

\begin{eqnarray}\myownlab{flcc}
X_{n-1}Y_{n-1}^mX_{n-1}&=&Y_n^mX_{n-1}^2=
Y_n^m+\delta Y_n^mX_{n-1}-\delta\lambda Y_n^me_{n-1}\\
&=&Y_n^m+\delta X_{n-1}Y_{n-1}^m-\delta\lambda Y_n^me_{n-1}\nonumber\\
&&\mbox{The last term can be  reduced using (\ref{P2a1})} \nonumber
\end{eqnarray}

The remaining cases (marked by $\ast$ in the table) are
\begin{eqnarray*}
Y_n^lY_{n-1}^mY_n^j &=& 
  X_{n-1}Y_{n-1}X^{-1}_{n-1}Y_{n-1}^mX_{n-1}Y_{n-1}^jX_{n-1}^{-1}\\
X_{n-1}Y_{n-1}^mY_n^j &=& 
  X_{n-1}Y_{n-1}^mX_{n-1}Y_{n-1}^jX_{n-1}^{-1}\\
Y_n^lY_{n-1}^mX_{n-1} &=& 
  X_{n-1}Y_{n-1}X^{-1}_{n-1}Y_{n-1}^mX_{n-1}
  \end{eqnarray*}
We note that we are dealing with sequences of generators where all
indices are equal. Hence we will suppress the index in further
calculations. Equations  (\ref{AKlem2eyx}) and (\ref{AKlem2xye})
imply that every such sequence containing $e$ is reducible
to $Y^teY^s$ and thus is of the standard form.
This motivates the following notation: 
We write  $a\sim b$ if 
$\exists c,\chi\;a=b+\chi c$ where $c$ contains $e$ 
and $\chi$ is some parameter. As a consequence 
the substitutions $X-\delta\leftrightarrow X^{-1}$  preserve this 
equivalence relation. 

To complete the proof it suffices to show that any finite 
sequence of the kind
$\cdots XY^{i_1}XY^{i_2}X\cdots$ is equivalent under $\sim$
to a sequence that contains at most two $X$ because if the sequence
contains none or only one $X$ it is in the standard form and if 
it contains exactly two $X$ it is either 
$XY^lXY^m\sim XY^lX^{-1}Y^m+\delta XY^{l+m}=
Y^l_{\cdot+1}Y^m+\delta XY^{l+m}$ 
or $Y^lXY^mX\sim	Y^lXY^mX^{-1}+\delta Y^lXY^{m}=
Y^lY^m_{\cdot+1}+\delta Y^lXY^{m}$.

The reducibility to sequences with at most two $X$ follows 
by induction from the following lemma:
There exists families of scalars $\alpha,\beta$ such that
\begin{equation}
XY^sXY^tX\sim\sum_{i,j}\alpha_{s,t}^{s,t} XY^iXY^j +
\sum_{i,j}\beta_{i,j}^{s,t} Y^iXY^j  
\label{AKredsim}
\end{equation}

We prove (\ref{AKredsim}) by induction on $s$. For $s=1$ we have
$XYXY^tX=Y^tXYX^2\sim Y^tXY-\delta Y^tXYX=Y^tXY-\delta XYXY^t$.
Assume that (\ref{AKredsim}) holds for $s$. We show it for $s+1$:
\begin{eqnarray*}
XY^{s+1}XY^tX&=&XYX^{-1}XY^sXY^t \\
&\sim&\sum_{i,j}\alpha^{s,t}_{i,j}  XYX^{-1}XY^iXY^j +
\sum_{i,j}\beta_{i,j}^{s,t} XYX^{-1}Y^iXY^j\\
&\sim&\sum_{i,j}\alpha^{s,t}_{i,j}  XY^{i+1}XY^j +
\sum_{i,j}\beta_{i,j}^{s,t} XYXY^iXY^j-\\&&
\delta \sum_{i,j}\beta_{i,j}^{s,t} XY^{i+1}XY^j
\end{eqnarray*}
The first and third summand are already in a form in which their 
contribution to $\alpha^{s+1,t}_{i+1,j}$ can be read off.
In the second summand we apply the induction hypothesis once again
\begin{eqnarray*}
\lefteqn{\sum_{i,j}\beta_{i,j}^{s,t} XYXY^iXY^j}\\
&\sim& \sum_{i,j}\beta_{i,j}^{s,t}\sum_{p,q}
\left(\alpha^{1,i}_{p,q} XY^pXY^{q+j}+\beta^{1,i}_{p,q}Y^pXY^{q+j} \right)
\end{eqnarray*}

We now establish the last statement of the proposition. 
Using the involution from remark \ref{AKinvolution} we see that we may
replace $X$ and $Y$ in $\Gamma_n$ by their inverses. Since $Y_n^{-1}$
is just a linear combination of powers of $Y_n$ we also may 
replace $Y_n$ by $Y_n^{-1}$ alone. Combining both operations 
replaces just $X_{n-1}$ by its 	inverse.
\end{bew}

The proposition implies  that $\AKBB_n$ is finite dimensional.

\begin{lemma} \label{lemeyes} There exist elements 
$R_{i,m}\in\AKBB_{i-1}$	such that $e_i{Y'}_i^me_i=R_{i,m}e_i$.
\begin{equation}
e_i{Y'}_i^lX_i=\kappa\lambda e_i{Y'}_i^{l-1}X_i{Y'}_i^{-1}
-\kappa\delta\lambda e_i{Y'}_i^{l-2}+
\kappa\delta\lambda R_{i,l-1}e_i{Y'}_i^{-1}	\myownlab{Reyx}
\end{equation}
\end{lemma}
\begin{bew}
To prove the first statement  one writes 
${Y'}_i^m=\sum_j a_jb_jc_j$ according to proposition \ref{wortsatz}
with $a_j,c_j\in\AKBB_{i-1},b_j\in\Gamma_i$. The claim is then obvious.

\begin{eqnarray}
e_i{Y'}_i^lX_i&=&e_i{Y'}_i^lX_iY_i'X_iX_i^{-1}{Y'}_i^{-1}=
e_iX_iY_i'X_i{Y'}_i^lX_i^{-1}{Y'}_i^{-1}\nonumber\\
&=&\kappa\lambda e_i{Y'}_i^{l-1}X_i^{-1}{Y'}_i^{-1}\nonumber\\  
&=&\kappa\lambda e_i{Y'}_i^{l-1}X_i{Y'}_i^{-1}-
\kappa\delta\lambda e_i{Y'}_i^{l-1}{Y'}_i^{-1}+
\kappa\delta\lambda e_i{Y'}_i^{l-1}e_i{Y'}_i^{-1}\nonumber\\
&=&\kappa\lambda e_i{Y'}_i^{l-1}X_i{Y'}_i^{-1}
-\kappa\delta\lambda e_i{Y'}_i^{l-2}+
\kappa\delta\lambda R_{i,l-1}e_i{Y'}_i^{-1}	
\end{eqnarray}
\end{bew}

\begin{satz} \myownlab{yswort}
In proposition \ref{wortsatz} one may replace $\Gamma_n$ 
by $\Gamma'_n:=\{1,e_{n-1},X_{n-1},{Y'}^{j}_n,j=1,\ldots,k-1\}$.
\end{satz}
\begin{bew}
We express an arbitrary element  $a$ in $\AKBB_n$
as $a=\sum_j f_jh_jg_j$ with $f_j,g_j\in\AKBB_{n-1},h_j\in\Gamma_n$.
We are finished if we can show that
$Y_n^i=\sum_s l_s^{(n)}\gamma_s^{(n)}r_s^{(n)}$ with
$\gamma_s^{(n)}\in\Gamma'_n, l_s^{(n)},r^{(n)}_s\in\AKBB_{n-1}$
since in this case we can simply substitute this expressions
for the $Y^n_i$ which appear among the $h_j$. 

We show $Y_n^i=\sum_s l_s^{(n)}\gamma_s^{(n)}r_s^{(n)}$ by induction.
The case $n=1$ is trivial. Now assume that the formula holds for
$n-1$. 
\begin{eqnarray*}
Y_n^i&=&X_{n-1}Y_{n-1}^iX_{n-1}^{-1}=
\sum_s X_{n-1}l_s^{(n-1)}\gamma_s^{(n-1)}r_s^{(n-1)}X_{n-1}^{-1}\\
&=&	\sum_s l_s^{(n-1)}X_{n-1}\gamma_s^{(n-1)}X_{n-1}^{-1} r_s^{(n-1)}
\end{eqnarray*}
The cases $\gamma_s^{(n-1)}\in\{1,e_{n-2},X_{n-2}\}$
are easily reduced using lemma \ref{lemma2}. It remains to
investigate the case $\gamma_s^{(n-1)}={Y'}_{n-1}^{j}$.
\begin{eqnarray}
X_{n-1}{Y'}_{n-1}^{j}X_{n-1}^{-1}&=&
X_{n-1}{Y'}_{n-1}X_{n-1}X_{n-1}^{-1}{Y'}_{n-1}^{j-1}X_{n-1}^{-1}	
\nonumber\\
&=&{Y'}_n(X_{n-1}-\delta+\delta e_{n-1}){Y'}_{n-1}^{j-1}X_{n-1}^{-1}
\nonumber\\
&=&{Y'}_nX_{n-1}{Y'}_{n-1}^{j-1}X_{n-1}^{-1}-
\delta {Y'}_n{Y'}_{n-1}^{j-1}X_{n-1}^{-1}+\nonumber\\&&
\delta {Y'}_ne_{n-1}{Y'}_{n-1}^{j-1}X_{n-1}^{-1}	 \myownlab{props2}
\end{eqnarray}
The second summand is
$-\delta {Y'}_{n-1}^{j-1}X_{n-1}{Y'}_{n-1}$ which is already of the
standard form. 
The third summand is 
\begin{eqnarray*}
\lefteqn{\delta Y'_ne_{n-1}{Y'}_{n-1}^{j-1}X_{n-1}^{-1}
=\delta\lambda X_{n-1}Y'_{n-1}e_{n-1}{Y'}_{n-1}^{j-1}
(X_{n-1}-\delta+\delta e_{n-1})}\\
&=&\delta\lambda\kappa {Y'}_{n-1}^{-1}e_{n-1}{Y'}_{n-1}^{j-1}X_{n-1}
-\delta^2\lambda\kappa {Y'}_{n-1}^{-1}e_{n-1}{Y'}_{n-1}^{j-1}
+ \\&&\delta^2\lambda\kappa {Y'}_{n-1}^{-1}e_{n-1}{Y'}_{n-1}^{j-1}e_{n-1}
\end{eqnarray*}
Here the last summand is reduced using the formula
for $e_i{Y'}_i^me_i$ from lemma \ref{lemeyes} 
while the  first summand needs
(\ref{Reyx}). The middle summand is already of the standard form. 

The first summand of (\ref{props2}) 
is reduced by iteration.
\end{bew}

We continue our study of words in $\AKBB_n$ by cutting down the size
of sets that linearly generate the algebra.

\begin{lemma}		\label{ysystem}
$\AKBB_n$ is linearly spanned by the set $S_n$ which is recursively
defined:
\begin{eqnarray*}
S_1&:=&\{Y^i\mid i=0,\ldots,k-1\}\\
S_n&:=&\Gamma_1\cdots\Gamma_nS_{n-1}
\end{eqnarray*}
It suffices to take out of  $\Gamma_1\cdots\Gamma_n$ those elements
that are of the following form:
\[ Y^{m_1}_{l_1}\cdots Y^{m_s}_{l_s}X_i\cdots X_je_{j+1}\cdots e_n,\qquad
m_t\in\{0,\ldots,k-1\},m_s=i,l_1<\cdots<l_s\]
Here we have  $1\leq i\leq n$ and $i-1\leq j\leq n$ so that
the chains of  $X$ and $e$	may be empty.
\end{lemma}
\begin{bew}
Proposition \ref{wortsatz} yields the following
representation of $\AKBB_n$:
\begin{eqnarray*}
\AKBB_n&=&{\rm span}\AKBB_{n-1}\Gamma_n\AKBB_{n-1}
={\rm span}\AKBB_{n-2}\Gamma_{n-1}\AKBB_{n-2}\Gamma_n\AKBB_{n-1}\\
&=&{\rm span}\AKBB_{n-2}\Gamma_{n-1}\Gamma_n\AKBB_{n-1}\\
&=&{\rm span}\Gamma_1\cdots\Gamma_n\AKBB_{n-1}
\end{eqnarray*}
To establish the second statement we consider the
$Y_j^m$ that appears at the leftmost position in a chain
$Z_i\cdots Z_{j-1}Y_j^mZ_{j+1}\cdots Z_n$
of generators $Z_s\in\Gamma_s$. 
Then $Z_i\cdots Z_{j-1}$ consists only of $e$ and $X$ 
and hence it can be commuted to the right and be absorbed 
in  $\AKBB_{n-1}$. Similarly  $e$ and $X$ that appear between
two  $Y_\cdot$ can be commuted to the right.
Iterating this argument we obtain only chains of the form
$Y^{m_1}_{i_1}\cdots Y^{m_s}_{i_s}Z_{j+1}\cdots Z_n,i_1<\cdots<i_s$.

If   $e_iX_{i+1}$ appears in such a chain
it may be converted to $e_iX_{i+1}=e_ie_{i+1}X_i^{-1}$. 
The $X_i^{-1}$ can the be absorbed in $\AKBB_{n-1}$.
Hence all $X$ have to appear to the left of all $e$.
\end{bew}

A similar proof establishes a related lemma
using the $Y'_i$:

\begin{lemma}	\label{yssystem}
$\AKBB_n$ is linearly spanned by $S'_n$: 
\begin{eqnarray*}
S'_1&:=&\{{Y'}^i\mid i=0,\ldots,k-1\}\\
S'_n&:=&\Gamma'_1\cdots\Gamma'_nS_{n-1}
\end{eqnarray*}
From $\Gamma'_1\cdots\Gamma'_n$ only elements of the following form
are needed:
\[ {Y'}^m_iX_i\cdots X_je_{j+1}\cdots e_n,\qquad
m=0,\ldots,k-1\]
The chains of $x$ and $e$ may be empty. \end{lemma}

\section{Graphical interpretation and classical limit}
\myownlab{clsec}

The very definition of  $\AKBB_n$ is motivated by knot theory
as was vaguely explained in section \ref{AKdefsec}. 
Here we fill in the details.

Consider the free  $R$ algebra
($R$ may denote any commutative ring)
of isotopy classes of ribbons in $(\RR^2-\{0\})\times [0,1]$
where $n$ ribbons end at the upper and lower  plane each.
The ribbons touch these planes in small intervals
which have as their lower starting point one of the set 
$\{1,2,\ldots,n\}\times0\times\{0,1\}$.  Closed components are allowed.
Multiplications is given by putting graphs on top of each other.
This forms the algebra of cylinder tangles.

The pictures on the right hand side of figure \ref{generat}
may now be easily interpreted as regular diagrams of 
such cylinder tangles. We need to specify the total number of strands
in these pictures. Thus, we write  $X^{(G)}_{i,n}$,
$e^{(G)}_{i,n}$ and $Y^{(G)}_{i,n}$ for the generators that act
at the $i$-th of $n$ strands. 
Let ${G\AKBB_n}'(R)$ be the sub-algebra of the algebra of cylinder tangles
that is generated by  $X_{i,n}^{(G)},e_{i,n}^{(G)},Y_{1,n}^{(G)},
1\leq i\leq n-1$. Each isotopy class thus has a representative 
that is a product in these generators.
We define $G\AKBB(R)$ 
(where $R$ is now as in   definition \ref{definAK}) 
to be the quotient of this algebra by 
skein relations that result from
(\ref{def4}), (\ref{def5}) and (\ref{AKdef7}) by 
replacing  $X_i, e_i, Y$ by 
$X_{i,j}^{(G)},e_{i,j}^{(G)},Y_{1,j}^{(G)}$. 
Here, we don't restrict  $j$ so that it may be
greater than $n$. This is necessary to account for the fact 
that by introducing maxima and minima the number of strands that
intersect some horizontal plane may be arbitrary.
The remaining relations of $\AKBB_n(R)$ have obvious 
topological content so that we have a surjective
morphism \label{psimorph}
$\Psi_n:\AKBB_n(R)\rightarrow G\AKBB_n(R)$. 
We remark that this graphical algebra is not defined in terms of a basis
but in terms of generators and relations. However, some of the 
relations are not stated explicitly.
The existence of $\Psi_n$, however, shows that
 the statements of section \ref{wordsec} carry over.
However, we have to keep in mind the possibility that
 $\Psi_n$ could fail to be injective.

The graphical interpretation suggests 
special settings for $\kappa$.
Recall that $\lambda$ amounts to a twist of the ribbon.
If we interpret $Y^{(G)}$ as a ribbon band that lies flat in the projection 
plane then  we should have $\lambda=\kappa$. On the other hand, if
the transversal vector field of the ribbon is always oriented towards the 
cylinder axes we should have $\kappa=1$. However, we can (and will)
decide to keep $\kappa$ free by renormalising $Y$.

The classical limit of tangle algebra is a specialisation in which
braidings degenerate to permutations. 
We define  ${\rm BP}^k_n(R)$ in its own right as 
algebra of  Brauer graphs \cite{we1} where each arc carries
an element of  $\ZZ_k$. We visualise this as 
dotted Brauer graphs, i.e. ${\rm BP}^k_n(R)$ is the free
 $R$ module of dimension  $k^n(2n-1)!!$
that has as basis the set of Brauer graphs 
where each arc carries at most  $k-1$ points.
We require that vertical arcs have no extrema with respect to 
the height function and that horizontal arcs have exactly one extremum.
Furthermore, we demand that the dots of vertical arcs are
concentrated at the left endpoint.

Multiplication is given as for graphs.
Dots may flow along an arc and may cross another arc. 
If a dot traverses an extremum it gets replaced by $k-1$ dots. 
Dot numbers are reduced modulo $k$.
Using this we may isolate cycles and concentrate  dots on their leftmost
position.  Such a cycle with $i$ dots on it may be deleted at 
the expense of a factor $A_i$. Dots on vertical arcs may be brought 
to the lower endpoint and thereafter the arc may be straightened. 
Similarly, dots on horizontal arcs may be concentrated according to our convention.
Just as in the case of ordinary  Brauer graphs 
we see that ${\rm BP}^k_n(R)$ is generated by
$X_{i,n}^{(G)},e_{i,n}^{(G)},Y_{1,n}^{(G)}$ (where  $X^{(G)}_{i,n}$
is to be understood as a permutation two-cycle).

Lets compare ${\rm BP}^k_n$ with the classical limit of
 $\AKBB_n(R_1)$.
\begin{de} The classical limit of  $\AKBB_n(R_1)$ is
defined to be the algebra
\begin{eqnarray}
C\AKBB_n&:=&\AKBB_n(R_1)\otimes_{R_1}(R_1/J_c)\\
J_c&:=&	\label{jcdef} (\kappa-1,\lambda-1,\delta,q_0-1,q_1,\ldots,q_{k-1})
\subset R_1
\end{eqnarray}
The new ground ring  $R_1/J_c$ is denoted by  $R_c$.
\end{de}
Note that $(\kappa-1,\lambda-1,q_0-1,q_1,\ldots,q_{k-1})$
viewed as ideal in $R_0$ contains the consistency ideal $c$
because in the limit $Y^{-1}=Y^{k-1},
\overline{q}_{k-1}=1,\overline{q}_i=0$
and hence (\ref{consisteq})  becomes trivial.
Thus,  $R_c$ is the quotient of  $R_0$ by this ideal.

In  $C\AKBB_n$ we have  $X_i=X_i^{-1}$ and hence  
$Y_i^j=Y_i^{(j)}={Y'}_i^{j}$. 
An important consequence is that 
$Y'_i$ behaves natural with respect to the braidings $X_i$. 
In the system  $S'_n$ from 
lemma \ref{yssystem} we may read  $Y'$ as $Y$.
Using this we are going to prove that
$\AKBB_n$ is linearly spanned by a set of  elements of the form
$\alpha\beta\gamma$, where  $\alpha$ is a product of
 $Y_\cdot$,  $\gamma$ is a product of  $Y_\cdot^{-1}$ and 
$\beta$ is an element of a basis of the A-type BMW algebra 
$\BA_n$ .
The proof is by induction on $n$, so assume the claim is already 
shown for $n-1$.
It suffices to show that all $Y_i$ which appear on the left
of the generating system  $S_{n-1}$ of $\AKBB_{n-1}$ can be moved to 
the left through the left chain or that it can (in negated form)
be moved to right of 	  $\AKBB_{n-1}$.
We investigate the various arising cases.
In the first case    $e_{n-1}Y_{n-1}$ appears.
We rewrite it according to
\begin{eqnarray*}
e_{n-1}Y_{n-1}&=&
e_{n-1}Y_{n-1}X_{n-1}Y_{n-1}{Y_{n-1}}^{-1}X_{n-1}^{-1}
= e_{n-1}{Y_{n-1}}^{-1}X_{n-1}^{-1}\\
&=& e_{n-1}X_{n-1}{Y_{n-1}}^{-1}X_{n-1}^{-1}
= e_{n-1}{Y_n}^{-1}
\end{eqnarray*}
The $Y_n^{-1}$ may then be moved to the right.
If $e_ie_{i+1}Y_i=e_iY_ie_{i+1}$	occurs 
a twofold application of this result shows that
$ Y_{i+2}$ may be moved to the left.
The only remaining situation 
$X_iY_i=Y_{i+1}X_i$ is trivial.

In each 
step of the recursive construction of  $S'_n$ 
only one  additional $Y^m_i$ can occur and these occurences stick together
in the above process. The dimension of $C\AKBB$ 
is therefore at most $k^n$ times the dimension of
the ordinary Brauer algebra:  ${\rm dim}C\AKBB_n\leq k^n(2n-1)!!$. 

\begin{lemma} The algebras
$C\AKBB_n$ and ${\rm BP}^k_n(R_c)$ are isomorphic.
\end{lemma}
\begin{bew}
We define the morphism  
$\chi_n:C\AKBB_n\rightarrow {\rm BP}^k_n(R_c)$ that maps
$e_i\mapsto e_i^{(G)},X_i\mapsto X_i^{(G)}$ and
$Y$ to a dot on the frist strand. It is easy to see that
this is a morphism
(It is relation (\ref{AKdef8}) that requires the somewhat
strange minmum/maximum rule.). It is surjective. Injectivity 
may be seen by looking at the dimension of these algebras.
\end{bew}

\begin{lemma}
\label{gheckelem}
The quotient of  $G\AKBB_n$ by the ideal $I^{(G)}$ generated by  
$e_{1,n}^{(G)}$ is isomorphic to the  Ariki-Koike Algebra $\AK^k_n$.
\end{lemma}
\begin{bew}
A graph is of the form  $ae_{1,n}^{(G)}b$
if and only if it contains horizontal arcs.
The quotient consists hence of those graphs that have only
vertical arcs. It is therefore the group of ribbon braids
in the cylinder. The relations of this group are known 
to be a subset of the relations of the Ariki-Koike algebra. 
The remaining relations follow from the imposed skein relations.
\end{bew}

At this point the importance of the index $n$ (total number of strands)
of the generator $e^{(G)}_{i,n}$ becomes obvious. 
Without fixing the total number of strings the ideal
would be the whole algebra because minima and maxima can be introduced
within the isotopy class of any diagram (cf. figure \ref{inklu} right). 
On the other hand we have avoided to restrict the 
number of strands when defining the skein relations.
Using this we obtain the following lemma.

\begin{lemma} The map 
$G\AKBB_n(R)\rightarrow G\AKBB_{n+2}(R),a\mapsto x^{-1}ae_{n+1}$
is injective for $n\geq1$. \myownlab{gwenzlincl}
\end{lemma}
\begin{bew}
By deforming the  $n$-th strand of a graph $a$ 
we may generate maxima and minima
as shown in figure \ref{inklu} (on the left).
Thus, locally, we obtain   $ae_{n+1}$. If $a$ is in the kernel then
this vanishes and hence $a=0$.
\end{bew}

   \unitlength1mm
 \begin{figure}[ht]
\begin{picture}(120,30)

\linethickness{0.1mm}

\put(10,10){\line(1,0){10}}
\put(10,10){\line(0,1){10}}
\put(20,10){\line(0,1){10}}
\put(10,20){\line(1,0){10}}
\put(12,12){\mbox{$a$}}

\put(11,20){\line(0,1){8}}
\put(13,20){\line(0,1){8}}
\put(15,20){\line(0,1){8}}
\put(17,20){\line(0,1){8}}
\put(19,20){\line(0,1){8}}

\put(11,2){\line(0,1){8}}
\put(13,2){\line(0,1){8}}
\put(15,2){\line(0,1){8}}
\put(17,2){\line(0,1){8}}
\put(19,2){\line(0,1){8}}

\put(30,12){\mbox{$=$}}

\put(40,10){\line(1,0){10}}
\put(40,10){\line(0,1){10}}
\put(50,10){\line(0,1){10}}
\put(40,20){\line(1,0){10}}
\put(42,12){\mbox{$a$}}

\put(41,20){\line(0,1){8}}
\put(43,20){\line(0,1){8}}
\put(45,20){\line(0,1){8}}
\put(47,20){\line(0,1){8}}
\put(49,20){\line(0,1){4}}

\put(41,2){\line(0,1){8}}
\put(43,2){\line(0,1){8}}
\put(45,2){\line(0,1){8}}
\put(47,2){\line(0,1){8}}
\put(49,6){\line(0,1){4}}

\put(52,6){\oval(6,6)[b]}
\put(52,24){\oval(6,6)[t]}

\put(55,10){\line(0,-1){4}}
\put(55,20){\line(0,1){4}}

\put(58,10){\oval(6,6)[t]}
\put(58,20){\oval(6,6)[b]}

\put(61,10){\line(0,-1){10}}
\put(61,20){\line(0,1){10}}


\put(90,5){\line(0,1){20}}
\put(94,15){\mbox{$=$}}

\put(100,5){\line(0,1){5}}
\put(103,10){\oval(6,6)[t]}
\put(109,10){\oval(6,6)[b]}

\put(112,10){\line(0,1){10}}

\put(100,25){\line(0,-1){5}}
\put(103,20){\oval(6,6)[b]}
\put(109,20){\oval(6,6)[t]}

\end{picture}
\caption{\myownlab{inklu} }
\end{figure}

\section{Conditional expectation and Markov trace}
\myownlab{trsec}

The graphical calculus as well as the relationship with the
A-type BMW algebra suggest that there should exist 
a Markov trace on $\AKBB_n$. 
We follow Wenzl's original approach \cite{we2}
as close as possible.

The constructions of this section can equally well be carried out
for $\AKBB_n$  and for its graphical counterpart $G\AKBB_n$.
Notationally, however, we'll stick to the former case.

The fundamental hypothesis for the following construction
is:
\begin{hypo}
The map  $\AKBB_n\rightarrow\AKBB_{n+2},a\mapsto x^{-1}ae_{n+1}$
is injective. \myownlab{wenzlincl}
\end{hypo}
Lemma \ref{gwenzlincl} has shown that this hypothesis is valid for
the graphical algebra.

Let $w=w_1\gamma w_2\in\AKBB_{n+1}$ where 
$w_i\in\AKBB_n,\gamma\in\Gamma_{n+1}$. 
Then we have
$e_{n+1}we_{n+1}=w_1e_{n+1}\gamma e_{n+1}w_2=sw_1w_2e_{n+1}$
with a factor  $s$ that assumes the values
$s=x,1,\lambda^{-1},A_m$ if  $\gamma=1,e_n,X_n,Y_{n+1}^m$.
Hypothesis \ref{wenzlincl} guarantees that the following map
is well defined.
\begin{de}						   Let
$\epsilon_n:\AKBB_{n+1}\rightarrow\AKBB_n$ 
be defined by
$e_{n+1}ae_{n+1}=:x\epsilon_n(a)e_{n+1}$.
\end{de}

Obviously, we have $\epsilon_n(w_1aw_2)=w_1\epsilon_n(a)w_2$ 
if $w_i\in\AKBB_n$.
Moreover,   (\ref{lem1l}) implies
$e_{n+1}=e_{n+1}e_ne_{n+1}=x\epsilon_n(e_n)e_{n+1}$ 
and thus $\epsilon_n(e_n)=x^{-1}$. 
Similarly,  (\ref{def5}) implies
$e_{n+1}=\lambda^\pm e_{n+1}X^\pm_ne_{n+1}=
\lambda^\pm x\epsilon_n(X^\pm_n)e_{n+1}$ and thus 
$\epsilon_n(X_n^\pm)=x^{-1}\lambda^\mp$.
From relation   (\ref{AKlem2ds})  we deduce
$e_{n+1}=A_m^{-1}e_{n+1}Y_{n+1}^me_{n+1}=
A_m^{-1}x\epsilon_n(Y_{n+1})e_{n+1}$ and hence
$\epsilon_n(Y_{n+1}^m)=A_mx^{-1}$. 

Iterating the conditional expectation yields 
a map which will turn out to be a Markov trace.

\begin{lede}\label{tracedef}
The iterated application of the conditional expectation
is denoted by $\tr(a):=\tr(\epsilon_{n-1}(a)), tr(1):=1$
and fulfils
$\tr(e_n)=\epsilon_n(e_n)=x^{-1},\quad
              \tr(X_n^\pm)=\epsilon_n(X_n^\pm)=x^{-1}\lambda^\mp,\quad
              \tr(Y_{n+1}^m)=\epsilon_n(Y^m_{n+1})=A_mx^{-1}$
  
\end{lede}

\begin{lemma} \myownlab{trlem1}
For any $w_1,w_2\in\AKBB_n,\gamma\in\Gamma_{n+1}$  we have
$\tr(w_1\gamma w_2)=\tr(\gamma)\tr(w_1w_2)$ and
$\epsilon_n(w_1\gamma w_2)=\tr(\gamma)ab$.
\end{lemma}
\begin{bew} 
The first statement follows from the second which is shown by the 
following calculation.
$x\epsilon_n(w_1\gamma w_2)e_{n+1}=e_{n+1}w_1\gamma w_2e_{n+1}=
w_1e_{n+1}\gamma e_{n+1}w_2=w_1x\epsilon_n(\gamma)e_{n+1}w_2=
w_1w_2x\epsilon_n(\gamma)e_{n+1}$. 
\end{bew}

\begin{lemma}
$\forall a\in\AKBB_n\quad
\epsilon_n(X_n^{-1}aX_n)=\epsilon_n(X_naX_n^{-1})=\epsilon_n(e_nae_n)=\epsilon_{n-1}(a)$
\end{lemma}
\begin{bew}
Let $a=w_1\gamma w_2\in\AKBB_n,w_i\in\AKBB_{n-1},\gamma\in\Gamma_n$.
Multiplying by $xe_{n+1}$  we obtain:
\begin{eqnarray*}
\lefteqn{x\epsilon_n(X_n^{-1}w_1\gamma w_2X_n)e_{n+1}=
x\epsilon_n(X_nw_1\gamma w_2X_n^{-1})e_{n+1}=}\\
&=&x\epsilon_n(e_nw_1\gamma w_2e_n)e_{n+1}=
x\epsilon_{n-1}(w_1\gamma w_2)e_{n+1}
\end{eqnarray*}
Omitting the arbitrary factors  $w_1,w_2$ yields:

$e_{n+1}(X_n^{-1}\gamma X_n)e_{n+1}=e_{n+1}(X_n\gamma X_n^{-1})e_{n+1}=
e_{n+1}(e_n\gamma e_n)e_{n+1}=x\tr(\gamma)e_{n+1}$.

This is checked by analysing the cases for the various values of $\gamma$
successively. For $\gamma=1$ nothing is to be shown.
For $\gamma=e_{n-1}$ we have
\begin{eqnarray*}
\lefteqn{e_{n+1}(X_n^{-1}e_{n-1} X_n)e_{n+1}=e_{n+1}(X_ne_{n-1} X_n^{-1})e_{n+1}=
e_{n+1}(e_ne_{n-1} e_n)e_{n+1}=xx^{-1}e_{n+1}}\\
&\Leftrightarrow&
e_{n+1}(X_{n-1}e_{n} X_{n-1}^{-1})e_{n+1}=e_{n+1}(X_{n-1}^{-1}e_{n} X_{n-1})e_{n+1}=
e_{n+1}e_ne_{n+1}=e_{n+1}\end{eqnarray*}

The case $\gamma=Y_n^m$ yields
\begin{eqnarray*}
\lefteqn{e_{n+1}(X_n^{-1}Y_n^m X_n)e_{n+1}=e_{n+1}(X_nY^m_n X_n^{-1})e_{n+1}=
e_{n+1}(e_nY^m_n e_n)e_{n+1}=x\tr(Y^m_n)e_{n+1}}\\
&\Leftrightarrow
e_{n+1}(X_n^{-1}Y^m_n X_n)e_{n+1}=e_{n+1}Y^m_{n+1}e_{n+1}=
e_{n+1}(e_nY^m_n e_n)e_{n+1}=A_me_{n+1}
\end{eqnarray*}
We rewrite the first expression to obtain:

$e_{n+1}X_n^{-1}Y^m_nX_ne_{n+1}=e_{n+1}e_nX_{n+1}Y^m_nX_ne_{n+1}=
e_{n+1}e_nY^m_nX_{n+1}X_ne_{n+1}=
e_{n+1}e_nY^m_ne_nX_{n+1}X_n=A_me_{n+1}e_nX_{n+1}X_n=A_me_{n+1}$.

The last case is $\gamma=X_{n-1}$. 
\begin{eqnarray*}
\lefteqn{e_{n+1}(X_n^{-1}X_{n-1} X_n)e_{n+1}=
e_{n+1}(X_nX_{n-1} X_n^{-1})e_{n+1}=}\\
&&=e_{n+1}(e_nX_{n-1} e_n)e_{n+1}=x\tr(X_{n-1})e_{n+1}\\
&\Leftrightarrow&
e_{n+1}(X_{n-1}X_{n} X_{n-1}^{-1})e_{n+1}=
e_{n+1}(X_{n-1}^{-1}X_{n} X_{n-1})e_{n+1}=\\&&\quad=
e_{n+1}(\lambda^{-1}e_n)e_{n+1}=\lambda^{-1}e_{n+1}\\
&\Leftrightarrow&
X_{n-1}e_{n+1}X_ne_{n+1}X_{n-1}^{-1}=X_{n-1}^{-1}e_{n+1}X_ne_{n+1}X_{n-1}=
\lambda^{-1}e_{n+1}=\lambda^{-1}e_{n+1}\\
&\Leftrightarrow&
X_{n-1}\lambda^{-1}e_{n+1} X_{n-1}^{-1}=
X_{n-1}^{-1}\lambda^{-1}e_{n+1}X_{n-1}=
\lambda^{-1}e_{n+1}
\end{eqnarray*}
\end{bew}

Just as in \cite{we2} we have the trace property in the semi-simple case.

\begin{lemma}\myownlab{wenzlspur} 
If   $I_{n+1}$  is semi-simple and  $\tr$ is a trace on  $\AKBB_n$
then $\tr$ is  a trace on  $\AKBB_{n+1}$. \myownlab{wspurlem}
\end{lemma}
\begin{bew}
It suffices to show that
$\tr(uv)=\tr(vu)\forall u,v\in\AKBB_{n+1}$. 
If one of the factors (say $u$) is contained in  $\AKBB_n$
this is easily seen:
$\tr(uv)=\tr(\epsilon_n(uv))=\tr(u\epsilon_n(v))=
\tr(\epsilon_n(v)u)=\tr(\epsilon_n(vu))=\tr(vu)$.

According to proposition  \ref{wortsatz} 
we may write  $u,v\in\AKBB_{n+1}$ in the form
\begin{eqnarray}
u&=&u_1+u_2Y'_{n+1}+u_3e_nu_4+u_5X_nu_6\\
v&=&v_1+v_2Y'_{n+1}+v_3e_nv_4+v_5X_n^{-1}v_6
\end{eqnarray}
Since  $\tr$ is linear it suffices to investigate all
possible combinations of summands.
The calculations are similar to those in \cite{we2} and
\cite{rhobmw}. Thus we only give calculations for the cases that
involve $Y$.
Our first case is: 
$a=a_1Y'_{n+1},b=b_1Y'_{n+1}$.
\begin{eqnarray*}
\tr(ab)&=&\tr(a_1Y'_{n+1}b_1Y'_{n+1})=\tr(a_1{Y'}_{n+1}^2b_1)\\
&=&\tr({Y'}^2_{n+1}b_1a_1)=\tr(b_1Y'_{n+1}a_1Y'_{n+1})=\tr(ba)
\end{eqnarray*}

Next we look at
 $a=a_1e_na_2,b=a_3Y'_{n+1}$.
\begin{eqnarray*}
\tr(ab)&=&\tr(a_1\epsilon_n(e_na_2a_3Y'_{n+1}))=
  \tr(a_1\epsilon_n(e_nY'_{n+1})a_2a_3)\\
&=&\lambda\tr(a_1\epsilon_n(e_n{Y'_n}^{-1})a_2a_3)=
\lambda\tr(a_1\epsilon_n(e_n){Y'_n}^{-1}a_2a_3)\\
&=&\lambda x^{-1}\tr(a_1{Y'_n}^{-1}a_2a_3)=
   \lambda x^{-1}\tr(a_2a_3a_1{Y'_n}^{-1})\\
&=&\lambda\tr(a_2a_3a_1\epsilon_n(e_n){Y'_n}^{-1})=
\tr(a_2a_3a_1\epsilon_n(Y'_{n+1}e_n))\\
&=&\tr(\epsilon_n(a_3Y'_{n+1}a_1e_n)a_2)=\tr(ba)
\end{eqnarray*}
The case  $b=a_1e_na_2,a=a_3Y'_{n+1}$ is treated similarly.

Next case: $a=a_1X_na_2,b=b_1Y'_{n+1}$.
\begin{eqnarray*}
\tr(ab)&=&\tr(a_1\epsilon_n(X_nY'_{n+1})a_2b_1)=
\tr(a_1\epsilon_n(X_n^2Y'_{n}X_n)a_2b_1)\\
&=&\tr(a_1\epsilon_n(Y'_nX_n)a_2b_1)+
    \delta\tr(a_1\epsilon_n(Y'_{n+1})a_2b_1)
	-\delta\lambda\tr(a_1\epsilon_n(e_nY'_nX_n)a_2b_1)\\
&=&\tr(a_1Y'_nx^{-1}\lambda^{-1}a_2b_1)+
    \delta\tr(a_1\epsilon_n(Y'_{n+1})a_2b_1)
	-\delta\tr(a_1\epsilon_n(e_nY'_{n+1})a_2b_1)\\
&=&\tr(a_1Y'_nx^{-1}\lambda^{-1}a_2b_1)+
    \delta\tr(a_1\epsilon_n(Y'_{n+1})a_2b_1)
	-\delta\tr(a_1\epsilon_n(Y'_{n+1}e_n)a_2b_1)\\
&=&\tr(a_1\epsilon_n(X_nY'_n)a_2b_1)+
    \delta\tr(a_1\epsilon_n(Y'_{n+1})a_2b_1)
	-\delta\lambda\tr(a_1\epsilon_n(X_nY'_ne_n)a_2b_1)\\
&=&\tr(a_1\epsilon_n(Y'_{n+1}X_n)a_2b_1)=
\tr(b_1\epsilon_n(Y'_{n+1}a_1X_n)a_2)=\tr(ba)
\end{eqnarray*}

The cases where
$e_n$ is matched with $X_n$ or $X^{-1}_n$ 
yield zero by semi-simplicity of the ideal. Namely, there is
a central idempotent  $z\in\AKBB_{n+1}$ such that
$z\AKBB_{n+1}\iso I_{n+1}$. Take $a\in I_{n+1}$ and thus $a=az$ and
$ab=azb=a(zb)$. Thus we may assume that  $b\in I_{n+1}$
as well.
But $a,b\in I_{n+1}$ implies that they are linear combinations of
terms of the form
$a=\sum_ia_ie_na'_i,b=\sum_ib_ie_nb'_i$ with
$a_i,a'_i,b_i,b'_i\in\AKBB_n$. 
Thus we are back in situations already treated.
\end{bew}

The trace  in the classical limit is $C\AKBB_n$ given 
by closing the strands from the right .
Hence, let $a\in{\rm BP}^k_n$ be a dotted Brauer graph and denote by
 $n_i(a)$ the number of cycles with $i$ dots on its closure.
 Then we have
 \begin{equation}
 \tr(a)=x^{-n}\prod_{i=0}^{k-1}A_i^{n_i(a)}
 \end{equation}

\begin{lemma}
The trace is nondegenerate on   $C\AKBB_n={\rm BP}^k_n$.
\end{lemma}
\begin{bew}
Let $\{v_i\mid i=1,\ldots,k^n(2n-1)!!\}$ be a linear basis
of dotted Brauer graphs.
It suffices to show 
${\rm det}(\tr(v_iv_j^\ast)_{i,j})\neq0$.

The involution $a\mapsto a^\star$ maps graphs to their
top-down mirror image and replaces each dot by $k-1$ dots.
Hence the closure of $aa^\ast$ is free of dots.
Now assume that $a$ has  $s$ upper (and hence $s$ lower)
horizontal arcs. Then there are $s$ cycles in $aa^\ast$.
Upon closing another $s$ cycles are produced from the 
the remaining horizontal arcs. The vertical arcs form a
permutation and $a^\ast$ contains the inverse permutation. 
Upon closing these $n-2s$ vertical arcs yield $n-2s$ 
cycles. The closure of   $aa^\ast$ has therfore a total of
 $n$ cycles and  $\tr(aa^\ast)=1$.

We now  specialise the ground ring:
$A_1:=\ldots :=A_{k-1}:=x^{-1}$. 
The trace is then a Laurent polynomial in $x$.
The choice for the $A_i$ implies that additional dots on an arc 
decrease the degree (in $x$) of the trace. 
If $\beta$ is an arc of   $a$ and  $b$ is any other graph
which does not contain an arc which is the mirror image of 
$\beta$. 
By considering the cases that $\beta$ is vertical and horizontal 
individually  one easily sees that the cycle in the closure of $ab$
which contains $\beta$ consists of more than two arcs from $a$ and $b$.
The closure of  $ab$ has therfore less cycles than
the closure of $aa^\ast$. We conclude that $b=a^\ast$
is the unique graph with highest $x$ degree of $\tr(ab)$.
We now consider the determinant of the trace.
\[ {\rm det}(\tr(v_iv_j^\ast)_{i,j})
=x^{-nk^n(2n-1)!!}{\rm det}
((x^{n_0(v_iv_j^\ast)}\prod_{s=1}^{k-1}x^{n_s(v_iv_j^\ast)/2})_{i,j})\] 
In each row the element at the diagonal is the unique element with
highest degree in $x$. Calculating the determinant thus yields 
a sum with a unique term of highest degree. Thus the 
determinant does not vanish.
\end{bew}

\section{The structure theorem}  \myownlab{mainsec}

In this section we determine the structure of  $\AKBB_n(K_1)$.
It will turn out to be semi-simple over this generic ground field.
We only need a few definitions on Young diagrams before we can state
the structure theorem.

A Young diagram $\lambda$ of size $n$ 
is a partition of the natural number $n$.
$\lambda=(\lambda_1,\ldots,\lambda_k),\sum_i\lambda_i=n,
\lambda_i\geq\lambda_{i+1}$. 
In the following we use ordered tuples of Young diagrams
$\underline{\lambda}=(\lambda^1,\ldots,\lambda^k)$ 
(cf. \cite{ariki}). The size of a tuple of Young diagrams is the sum of 
sizes of its components.
Let $\widehat{\Gamma}^k_n$ 
be the set of all $k$ tuples  of Young diagrams of sizes $n,n-2,\ldots$.

 \begin{satz} \myownlab{hauptsatz}
 \begin{enumerate}
\item  $\AKBB_n(K_1)$ is a semi-simple algebra isomorphic to 
$G\AKBB_n(K_1)$.
The simple components are indexed by
$\widehat{\Gamma}^k_n$. \myownlab{decomp}
\begin{equation}\AKBB_n=
\bigoplus_{\underline{\lambda}\in\widehat{\Gamma}^k_n}
 \AKBB_{n,\underline{\lambda}}\end{equation}
\item 
The Bratteli rule for restrictions of modules: A simple
$\AKBB_{n,\underline{\nu}}$ module 
$V_{\underline{\nu}},\underline{\nu}\in\widehat{\Gamma}^k_n$ 
decomposes into $\AKBB_{n-1}$ 
modules such that the $\AKBB_{n-1}$ module  
$\underline{\lambda}\in\widehat{\Gamma}^k_{n-1}$ occurs iff
$\underline{\lambda}$ may be obtained from $\underline{\nu}$ 
by adding or removing a box.
\item $\tr$ is a faithful trace. To every tuple of Young diagrams 
$\underline{\lambda}\in\widehat{\Gamma}^k_n$ 
there is an idempotent
$p_{\underline{\lambda}}$ and a non vanishing, rational function	 
$Q_{\underline{\lambda}}$
which does not depend on $n$ and satisfies
$\tr(p_{\underline{\lambda}})=Q_{\underline{\lambda}}/x^n$.
\end{enumerate}
\end{satz}

For the proof of the structure theorem we need  some  facts from
Jones-Wenzl theory of
inclusions of finite dimensional semi-simple algebras.

Let  $A\subset B\subset C$ 
be a unital embedding of finite dimensional semi-simple algebras and let
$\tr$  be a trace on $A, B$ 
that is compatible with the inclusion.
The associated conditional expectation is denoted by 
$\epsilon_A:B\rightarrow A,\tr(ab)=\tr(a\epsilon_A(b))$.
It is assumed that there is an idempotent $e\in C$ such that
$e^2=e,ebe=e\epsilon_A(b)\forall b\in B$ and  
$\varphi:A\rightarrow C,a\mapsto ae \mbox{ is injective}$.

Such a situation can be realized starting from an inclusion pair 
$A\subset B$ with a common faithful trace  $\tr$ 
and conditional expectation $\epsilon_A$.
We set 
$\widehat{C}:=\{\alpha:B\rightarrow B\mid \mbox{ linear},
\alpha(ba)=\alpha(b)a\forall a\in A,b\in B\}$.
The inclusion  $B\subset\widehat{C}$ 
 is given by
$b\mapsto\alpha_b,\alpha_b(b_1):=bb_1$. Here  $e$ is given by 
$e_A=\epsilon_A:B\rightarrow B$.
The sub-algebra of $\widehat{C}$ generated by  $B$ and  $e_A$ 
is denoted by $<B,e_A>$.
For this setup Wenzl has obtained the following results 
\cite[Theorem 1.1]{we2}
\begin{enumerate}
\item $<B,e_A>\iso{\rm End}_A(B)$\myownlab{towera}
\item The simple components of 
$A$ and $<B,e_A>$ are in 1-1 correspondence. The inclusion matrices of
$A\subset B\subset<B,e_A>$
are relatively transposed. 
If $p$ is a minimal idempotent in $A$ then $pe_A$ is a minimal idempotent
in $<B,e_A>$\myownlab{towerb}
\item $<B,e_A>\iso Be_AB$ \myownlab{towerc}
\item $<B,e>\iso<B,e_A>\oplus\widetilde{B}$ where $\widetilde{B}$
is a sub-algebra of $B$. \myownlab{towerd}
\item \ref{towerd} implies that the ideal generated by  $e$ in $C$ 
is isomorphic to $<B,e_A>$.
\myownlab{towere}
\end{enumerate}

We now  prove  the main theorem.

\begin{bew}
$\AKBB_0$ is simply the ground ring. 
Thus the proposition is true with
$\tr(p_{(\cdot,\ldots,\cdot)})=\tr(1)=Q_{(\cdot,\ldots,\cdot)}/x^0,
Q_{(\cdot,\ldots,\cdot)}=1$. 
The algebra $\AKBB_1$ is of dimension $k$ and has a basis 
$\{1,Y,\ldots Y^{k-1}\}$. 
It is commutative and semi-simple. The simple blocks are given by
the eigen spaces of $Y$. Existence of idempotents is clear.
The graphical version is isomorphic as a simple consequence of 
Turaev's result on the skein module of the solid torus \cite{turaevskein}.

We have to establish that the trace on 
 $\AKBB_1$ is nondegenerate. Now, suppose it were degenerate.
 Choose  $0\neq a=\sum_ma_mY^m\in\AKBB_1$ such that
  $\forall b\in\AKBB_1:\tr(ab)=0$. It is
$\tr(ab)e_1=x^{-1}e_1abe_1$. We have already noted
in the discussion  following definition \ref{defconsist} that the
ideal generated by  $e_1$ in $\AKBB_2$ is simple and that it is spanned by
$Y^ie_1Y^j$. The trace does not vanish on this module and hence it
is faithful. Furthermore the embedding $\AKBB_1\subset\AKBB_2$ 
is faithful and hence there must be an element  $b'$ 
in the ideal $I_2$ such that  $\tr(ab')\neq0$. 
Suppose $b'=\sum_{i,j}c_{i,j}Y^ie_1Y^j$.
We then have $x\epsilon_1(ab')e_2=\sum_{i,j,m}a_mc_{i,j}e_2Y^mY^ie_1Y^je_2
=\sum_{i,j,m}a_mc_{i,j}Y^mY^ie_2e_1e_2Y^j
=\sum_{i,j,m}a_mc_{i,j}Y^{m+i}e_2Y^j$ and hence 
$\tr(ab')e_1=x^{-2}e_1\sum_{i,j,m}a_mc_{i,j}Y^{m+i}Y^je_1=
x^{-2}\sum_{m,i,j}a_mc_{i,j}e_1Y^{m+i+j}e_1\neq0$. This is
the required contradiction.

Assume the proposition is shown by induction for $\AKBB_n$.

By induction assumption $\AKBB_n=G\AKBB_n$. 
We investigate the kernel of the inclusion
 $i:\AKBB_n\rightarrow\AKBB_{n+2}$ introduced in section \ref{trsec}.
Assume $i(a)=0$, then we have 
$0=\Psi_{n+2}(i(a))=i^{(G)}(a)$. Since  $i^{(G)}$
is injective according to lemma \ref{gwenzlincl}
it follows that $a=0$ and hence that $i$ is injective.
Thus we can use the results from section \ref{trsec} . 

We apply  Jones-Wenzl theory to the following situation:
$A=\AKBB_{n-1},B=\AKBB_n,C=\AKBB_{n+1},e=x^{-1}e_n,
\epsilon_A=\epsilon_{n-1}$.
This is possible because $A,B$ are semi-simple algebras with a
faithful trace by induction assumption.
All properties needed for $e$ have already been established.
Statement \ref{towera} of Jones-Wenzl theory asserts the semi-simplicity of
${\rm End}_A(B)\iso<B,e_A>$  
which is by \ref{towere} the ideal generated by $e$.
Thus $I_{n+1}$ is semi-simple.
The quotient algebra $\AKBB_{n+1}/I_{n+1}$ is the Ariki-Koike  algebra 
$\AK^{(k)}_{n+1}$ and is semi-simple according to \cite{ariki}.
Since we work over a field we can conclude (by looking at the radicals)
that $\AKBB_{n+1}$ is semi-simple and that it is isomorphic 
to the direct sum  
$\AKBB_{n+1}=I_{n+1}\oplus\AKBB_{n+1}/I_{n+1}$.
Statement \ref{towerb} asserts that the simple components of 
$I_{n+1}$ are indexed by $\widehat{\Gamma}^k_{n-1}$.
The simple components  of $\AK^{(k)}_{n+1}$ are indexed by 
tuples of Young diagrams of size $n+1$ (see \cite{ariki}). 

Consider the situation for the graphical algebra 
$G\AKBB_{n+1}$. By Jones-Wenzl theory we know that the ideals
 $I_{n+1}$ and $I_{n+1}^{(G)}$ are isomorphic. The quotient 
 $G\AKBB_{n+1}/I_{n+1}^{(G)}$ is by lemma \ref{gheckelem} 
isomorphic to the  Ariki-Koike algebra. Hence, we have
$G\AKBB_{n+1}=\AKBB_{n+1}$.
This completes the proof of point \ref{decomp} of the theorem.

The inclusion matrix for the part $I_{n+1}$ is the transpose
of the inclusion matrix of $\AKBB_{n-1}\subset\AKBB_n$. 
For the part $\AK^{(k)}_{n+1}$
the Bratteli rule follow from \cite{ariki}.

We have to show that $\tr$ is  faithful, i.e. that the
$Q$ functions don't vanish.
If $p_{\underline{\lambda}}\in\AKBB_{n-1}$ 
is a minimal idempotent in
$\AKBB_{n-1,\underline{\lambda}}$ then $x^{-1}p_{(\mu,\lambda)}e_n$ 
is a minimal idempotent in $\AKBB_{n+1}$ by point \ref{towerb}.
The trace of this idempotent is
$\tr(x^{-1}p_{\underline{\lambda}}e_n)=x^{-2}tr(p_{\underline{\lambda}})=
Q_{\underline{\lambda}}/x^{n-1+2}$.
Obviously, this is non vanishing (using the induction assumption). 
The idempotents of this kind are those of $I_{n+1}$. 
For the other idempotents (which are those of $\AKBB_{n+1}/I_{n+1}$)
the function $Q$ is defined by
$\tr(p_{\underline{\lambda}})=Q_{\underline{\lambda}}/x^n$.

To establish faithfulness of the trace we use the classical limit.
A minimal idempotent $p_{\underline{\lambda}}$ of $\AKBB_n$ yields
an idempotent in the classical limit described in section
\ref{clsec}. 
We know already that on the classical limit algebra the trace is
nondegenerate. Hence the function  $Q_{\underline{\lambda}}$	
does have a non vanishing classical limit and hence can't be zero
itself. 
\end{bew}

The Ariki-Koike-Birman-Wenzl algebra with $Y$ satisfying a quadratic
relation is of special interest and has been studied in \cite{rhobmw}. 

Naturally, one would like to study the algebra not only over 
the generic field but also with complex parameters.
The classical limit is then a point in parameter space and we know that
at this point the algebra is semi-simple. The functions $Q$  are apart from
a finite set of poles continuous and hence there is a neighbourhood
of  the classical point where the algebra 
is semi-simple. 
While some necessary conditions for semi-simplicity may be derived
easily from the knowledge of the Ariki-Koike algebra 
the determination of sufficient conditions has to await further studies.

   \unitlength1mm
 \begin{figure}[ht]
\begin{picture}(120,30)

\put(0,20){\mbox{${\rm BB}^2_0$}}
\put(0,10){\mbox{${\rm BB}^2_1$}}
 \put(0,0){\mbox{${\rm BB}^2_2$}}

\put(70,20){\mbox{$(\cdot,\cdot)$}}

\put(60,10){\mbox{$(\Box,\cdot)$}}
\put(80,10){\mbox{{\small $(\cdot,\Box)$}}}
\put(73,18){\line(1,-1){5}}
\put(73,18){\line(-1,-1){5}}

\put(65,0){\mbox{{\small $(\cdot,\cdot)$}}}
\put(75,0){\mbox{{\small $(\Box,\Box)$}}}
\put(35,0){\mbox{{\small $(\Box\Box,\cdot)$}}}
\put(50,0){\mbox{{\small $({\Box\atop\Box},\cdot)$}}}
\put(90,0){\mbox{{\small $(\cdot,{\Box\atop\Box})$}}}
\put(105,0){\mbox{{\small $(\cdot,\Box\Box)$}}}

\put(63,8){\line(-4,-1){20}}
\put(63,8){\line(1,-1){5}}
\put(63,8){\line(-1,-1){5}}
\put(63,8){\line(3,-1){15}}
\put(85,8){\line(-3,-1){15}}
\put(85,8){\line(1,-1){5}}
\put(85,8){\line(-1,-1){5}}
\put(85,8){\line(4,-1){20}}

\linethickness{0.2mm}

\end{picture}
\caption{\myownlab{brat} The first three lines of the
Bratteli digram of ${\rm BB}^2_n$}
\end{figure}
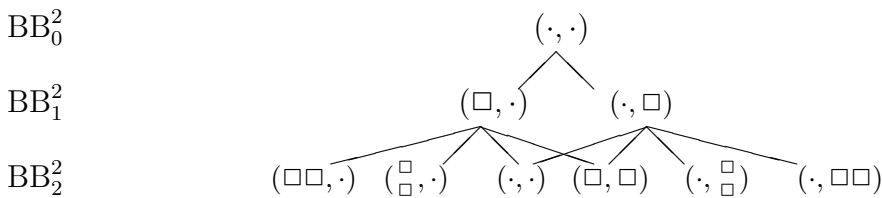

The Markov trace can be used to define a link invariant for 
links of B-type which are links in a solid torus. 
There is an analog of Markov's theorem for type B links found by S.
Lambrodopoulou in \cite{lamb}. It takes the same form as the 
usual Markov theorem,
i.e. two B-braids $\beta_1,\beta_2$ have isotopic closures
 $\hat{\beta}_1,\hat{\beta}_2$ if $\beta_1,\beta_2$
may transformed in one another by a finite sequence of moves of
the following two kinds: I Conjugation $\beta\sim\alpha\beta\alpha^{-1}$
and II $\alpha\sim\alpha\tau_n$ for $\alpha\in\ZB_n$.

This theorem implies that there exists an extension of the 
Kauffman polynomial 
to braids of B-type. Denote by $\pi:\ZB_n\rightarrow\AKBB_n$ the morphism
$\tau_i\mapsto X_i$. Then we obtain without any further proof
an invariant of the B-type link $\hat{\beta}$
that is the closure of a B-braid $\beta_\in\ZB_n$ by the following definition:
\begin{de}
The B-type Kauffman polynomial of a B-link $\hat{\beta}$ is defined to be
\begin{equation}
L(\hat{\beta},n):=x^{n-1}\lambda^{e(\beta)}\tr(\beta)\qquad\beta\in\ZB_n
\end{equation}
$e:\ZB_n\rightarrow\ZZ$ is the exponential sum with $e(X_i)=1,e(Y)=0$.
\end{de}

\end{document}